\documentclass[a4paper,12pt,oneside,reqno]{amsart}
\usepackage {amscd}
\usepackage {amsmath}
\usepackage {amssymb}
\usepackage [backend=biber,style=alphabetic,maxalphanames=4,sorting=nyt]{biblatex}
\usepackage {mathrsfs}
\usepackage {geometry}
\usepackage {setspace}
\usepackage {tikz}
\usepackage {tikz-cd}
\usetikzlibrary{arrows, decorations.pathmorphing}
\usepackage {float}
\usepackage {indentfirst}
\usepackage {enumitem}
\usepackage {amsthm}
\IfFormatAtLeastTF{2026-06-01}{}{\usepackage{aliascnt}}
\usepackage {hyperref}
\usepackage[capitalise]{cleveref}
\usepackage {mathabx}
\usepackage {diagbox}
\usepackage {graphicx}
\usepackage{longtable}
\usepackage{makecell}
\usepackage{etoolbox}
\usepackage{multirow}
\usepackage{adjustbox}
\preto\longtable{\setcounter{magicrownumbers}{0}}
\newcounter{magicrownumbers}

\graphicspath{ {./images/} }

\newcommand{\newsharedtheorem}[2]{%
  \IfFormatAtLeastTF{2026-06-01}{%
    \newtheorem{#1}[theorem]{#2}%
  }{%
    \newaliascnt{#1}{theorem}%
    \newtheorem{#1}[#1]{#2}%
    \aliascntresetthe{#1}%
  }%
}

\theoremstyle{plain}
\newtheorem{theorem}{Theorem}[section]
\newtheorem*{theorem*}{Theorem}
\newsharedtheorem{proposition}{Proposition}
\newsharedtheorem{corollary}{Corollary}
\newsharedtheorem{lemma}{Lemma}
\newsharedtheorem{conjecture}{Conjecture}

\theoremstyle{definition}
\newsharedtheorem{definition}{Definition}
\newsharedtheorem{proposition-definition}{Proposition-Definition}
\newsharedtheorem{question}{Question}
\newsharedtheorem{convention}{Convention}
\newsharedtheorem{condition}{Condition}

\theoremstyle{remark}
\newsharedtheorem{example}{Example}
\newsharedtheorem{construction}{Construction}
\newsharedtheorem{statement}{}
\newsharedtheorem{algorithm}{Algorithm}
\newtheorem*{remark}{Remark}

\crefname{theorem}{Theorem}{Theorems}
\crefname{proposition}{Proposition}{Propositions}
\crefname{corollary}{Corollary}{Corollaries}
\crefname{lemma}{Lemma}{Lemmas}
\crefname{conjecture}{Conjecture}{Conjectures}

\crefname{definition}{Definition}{Definitions}
\crefname{proposition-definition}{Proposition-Definition}{Proposition-Definitions}
\crefname{question}{Question}{Questions}
\crefname{convention}{Convention}{Conventions}
\crefname{condition}{Condition}{Conditions}

\crefname{example}{Example}{Examples}
\crefname{construction}{Construction}{Constructions}
\crefname{statement}{}{}
\crefname{algorithm}{Algorithm}{Algorithms}

\newcolumntype{C}{>{$}c<{$}} 

\title{Quantum periods and threefold flips}

\author{Yang He\textsuperscript{1}}
\email{heyang@bimsa.cn}
\author{Artan Sheshmani \textsuperscript{1,2,3}}
\email{artan@mit.edu}
\author{Shing-Tung Yau\textsuperscript{1,4}}
\address{\textsuperscript{1}Beijing Institute of Mathematical Sciences and Applications, No. 544, Hefangkou Village, Huaibei Town, Huairou District, Beijing 101408}
\address{\textsuperscript{2}Massachusetts Institute of Technology, IAiFi Institute, 77 Massachusetts Ave, 26-555. Cambridge, MA 02139}
\address{\textsuperscript{3}National Research University, Higher School of Economics, Russian Federation, Laboratory of Mirror Symmetry, NRU HSE, 6 Usacheva str., Moscow, Russia, 119048}
\address{\textsuperscript{4}Tsinghua University, No 30, Shuangqing road, Haidian district, Beijing}

\begin{document}

\begin{abstract}
    We study regularized quantum periods under birational contractions involving threefold flips. Let $g:X \dashrightarrow X'$ be a birational map between Fano threefolds with terminal quotient singularities which contracts exactly one prime divisor $E$ and extracts no divisor. We prove that
    $$
    \lim_{r \rightarrow +\infty}\widehat{G}_{X,rE}(t) =\widehat{G}_{X'}(t).
    $$
    The proof uses logarithmically smooth degenerations and positive discrepancies on good resolutions. Compatibility with orbifold invariants is proved using successive quotient-point extractions and the two-component orbifold degeneration formula.

    We also define stationary descendant invariants of projective terminal Fano threefolds by good resolutions and prove compatibility with a given projective $\mathbb{Q}$-smoothing. The period comparison extends to projective $\mathbb{Q}$-factorial terminal Fano endpoints under the same hypothesis on the birational map. Its proof compares the endpoints on a common resolution.
\end{abstract}

\maketitle

\section{Introduction}

    The quantum period of a Fano variety is a generating series of one-point genus-zero descendant Gromov-Witten invariants, organized by anticanonical degree. Understanding its behavior under birational transformations provides a connection between Gromov-Witten theory and the birational geometry of Fano varieties. This paper studies quantum periods of projective terminal Fano threefolds and their transformation under birational maps contracting a single prime divisor. It continues the study of quantum periods under threefold birational operations developed in \cite{Self_LG_Model,Self_LG_divisorial}.
    
    Two difficulties arise in this setting. First, a terminal threefold need not have quotient singularities, so the usual orbifold construction does not directly provide the required invariants. Second, a factorization of a birational map may pass through varieties that are not Fano, and individual flips need not preserve one-point Gromov-Witten invariants. Example \ref{example: GW change after flip} exhibits this failure explicitly. Our approach is to define the relevant invariants using resolutions and then compare the two endpoints on a common resolution.
    
    We first construct stationary descendant invariants for projective terminal Fano threefolds. Here the source curves are connected, there is at least one marking, and every insertion has positive even cohomological degree. The construction uses logarithmic Gromov-Witten invariants of a good resolution, summed over curve classes with prescribed pushforward and zero intersection number with every irreducible component of the exceptional divisor. We prove that these sums have finite support in each anticanonical degree and are independent of the resolution. Their one-point genus-zero specialization defines the quantum period. Compared with \cite{Self_LG_Model,Self_LG_divisorial}, we no longer require $\mathbb{Q}$-factoriality and the existence of a $\mathbb{Q}$-smoothing for this construction.
    
    We now state our main result.
    
\begin{theorem}[also \cref{thm:ordinary-period-limit}] 
    Let $X$ and $X'$ be projective $\mathbb{Q}$-factorial terminal Fano threefolds. Suppose that
    $$
    g:X\dashrightarrow X'
    $$
    is a birational map which contracts exactly one prime divisor $E$ and extracts no prime divisor. Then we have
    $$
    \lim_{r\rightarrow+\infty} \widehat{G}_{X,rE}^{\mathrm{res}}(t) = \widehat{G}_{X'}^{\mathrm{res}}(t).
    $$
    Here $\widehat{G}_{X,rE}^{\mathrm{res}}(t)$ is the regularized quantum period of $(X,rE)$ and $\widehat{G}_{X'}^{\mathrm{res}}(t)$ is the regularized quantum period of $X'$.
\end{theorem}

    For terminal quotient targets, Proposition \ref{prop:quotient-resolution-comparison} identifies the resolution invariants with the orbifold invariants with untwisted stationary insertions. \cref{thm:stationary-smoothing-compatibility} also proves compatibility with the nearby-fibre invariants of a given projective $\mathbb{Q}$-smoothing. Hence the resolution construction extends the quantum-period conventions used for quotient targets and for terminal targets equipped with a $\mathbb{Q}$-smoothing.
    
    The principal input is a virtual-dimension estimate for logarithmic relative invariants on smooth resolutions. Positive discrepancies, together with integral contact orders, force contributions with relative marked points to vanish for the stationary insertions under consideration. Applied to explicit degenerations of a smooth resolution and its blow-ups, this estimate gives the elementary comparisons needed for resolution independence. Boundary-compatible weak factorization then compares arbitrary good resolutions. The comparison with orbifold invariants is established separately, using successive quotient-point extractions and the two-component orbifold degeneration formula.

    The main theorem is proved by comparing the endpoints on a common resolution. The required control of the weighted coefficients is supplied by \cref{prop:ordinary-nonpositive-degree}.

\subsection{Organization of the paper}

    Section~2 recalls the orbifold Gromov-Witten invariants and quantum periods used in the comparison statements, together with the orbifold and logarithmic degeneration formulas.
    
    \cref{sec:scheme-comparisons} establishes the virtual-dimension and vanishing lemmas for degenerations.
    
    \cref{sec:stationary-terminal-invariants} defines the stationary invariants and constructs two log smooth degenerations used to prove resolution independence. It then establishes compatibility with quotient targets and with a given projective $\mathbb{Q}$-smoothing, and concludes with the definition of quantum periods for terminal Fano threefolds.
    
    \cref{sec:ordinary-terminal-case} begins with the toric flip example, develops the auxiliary relative comparisons and the nonpositive divisor-degree statement, and proves the main theorem.

\subsection*{Acknowledgment}

    The second author is supported by grants from Beijing Institute of Mathematical Sciences and Applications (BIMSA), the Beijing NSF BJNSF-IS24005, and the China National Science Foundation (NSFC) NSFC-RFIS program W2432008. He would also like to thank China's National Program of Overseas High Level Talent for generous support. Finally, he would like to thank NSF AI Institute for Artificial Intelligence and Fundamental Interactions at Massachusetts Institute of Technology (MIT) which is funded by the US NSF grant under Cooperative Agreement PHY-2019786.

\section{Preliminaries}

\subsection{The orbifold Gromov-Witten theory}\label{subsec:orbifold-gw}

    In this subsection, we recall some basic results about the orbifold Gromov-Witten theory, including the orbifold degeneration formula. The standard references in this topic are \cite{ChenRuan2001}, \cite{AGV08}, and \cite{AF16}.
    
\begin{convention}
    Throughout this article, a \emph{(smooth projective) orbifold} is a smooth separated Deligne-Mumford stack over $\mathbb{C}$ with a projective coarse moduli space. We say that an orbifold is Fano (resp. weak Fano) if the coarse moduli space is a Fano (resp. weak Fano) variety and is terminal if the coarse moduli space has terminal singularities.

    Let $X$ be a normal variety with only quotient singularities. We denote by $\mathfrak{X}$ its canonical covering stack. Thus $\mathfrak{X}$ is a smooth Deligne-Mumford stack with coarse moduli space $X$, and the morphism $\mathfrak{X} \rightarrow X$ is an isomorphism over the smooth locus of $X$. Whenever a variety with quotient singularities appears as a target of the orbifold Gromov-Witten theory, we implicitly replace it by its canonical covering stack. If $D$ is a Weil divisor on $X$, then the pullback of $D$ on $\mathfrak{X}$ is a Cartier divisor. We use the same letter $D$ for the pullback divisor on the canonical covering stack when there is no confusion.
\end{convention}

\begin{definition}[Inertia stacks]
    Let $\mathfrak{X}$ be a smooth Deligne-Mumford stack. Its inertia stack is
    $$
    \mathcal{I}\mathfrak{X} = \mathfrak{X} \times_{\mathfrak{X}\times \mathfrak{X}}\mathfrak{X}.
    $$
    Equivalently, an object of $\mathcal{I}\mathfrak{X}$ is a pair $(x,g)$, where $x \in \mathfrak{X}$ and \(g\in \mathrm{Aut}(x)\).  We also use the cyclotomic inertia stack
    $$
    \mathcal{I}_\mu(\mathfrak{X}) = \coprod_{r\ge 1}\underline{\mathrm{Hom}}^{\mathrm{rep}} (B\mu_r,\mathfrak{X}),
    $$
    and its rigidification $\overline{\mathcal{I}}_\mu\mathfrak{X}$.

    If $\Omega \subset \mathcal{I}_\mu(\mathfrak{X})$ is a connected component and \((x,g)\in \Omega\), write the eigenvalues of $g$ on $T_{x}\mathfrak{X}$ as
    $$
    \exp(2\pi i a_{1}),\cdots,\exp(2\pi i a_n), \qquad 0 \le a_{j}<1.
    $$
    The \emph{age} of $\Omega$ is
    $$
    \mathrm{Age}_{\mathfrak{X}}(\Omega)=\sum\limits_{j} a_{j}.
    $$
\end{definition}

\begin{theorem}[Reid--Tai criterion, {\cite[Theorem~4.11]{Reid87}}]\label{thm:reid-tai}
    Let $G\subseteq\operatorname{GL}_{n}(\mathbb{C})$ be a finite subgroup containing no quasi-reflections, that is, no nontrivial element fixes a hyperplane pointwise. For $g\in G$, write its eigenvalues as
    $$
    \exp(2\pi i a_1(g)),\ldots,\exp(2\pi i a_n(g)), \qquad 0\leq a_j(g)<1,
    $$
    and put
    $$
    \operatorname{Age}(g)=\sum\limits_{j=1}^{n}a_j(g).
    $$
    Then $\mathbb{C}^{n}/G$ has canonical, respectively terminal, singularities if and only if
    $$
    \operatorname{Age}(g)\geq1, \qquad\text{respectively}\qquad \operatorname{Age}(g)>1,
    $$
    for every $g\neq1$. In particular, if $X$ has terminal quotient singularities and $\mathfrak{X}$ is its canonical covering stack, then every nontrivial stabilizer element $\zeta$ satisfies
    $$
    \operatorname{Age}_{\mathfrak{X}}(\zeta)>1.
    $$
\end{theorem}

\begin{definition}[Twisted curves and twisted stable maps]
    A \emph{twisted curve} is a Deligne-Mumford stack $\mathfrak{C}$ whose coarse moduli space $C$ is a nodal marked curve, such that the stack structure is allowed only at marked points and nodes. Locally at a marked point, it has the form
    $$
        [\mathrm{Spec }\mathbb{C}[z]/\mu_r],
        \qquad
        \zeta\cdot z=\zeta z,
    $$
    and locally at a node, it has the balanced form
    $$
        [\mathrm{Spec }\mathbb{C}[x,y]/(xy)/\mu_r],
        \qquad
        \zeta\cdot(x,y)=(\zeta x,\zeta^{-1}y).
    $$
    Let $\mathfrak{X}$ be a smooth projective orbifold. A \emph{twisted stable map} to $\mathfrak{X}$ is a representable morphism
    $$
        f:\mathfrak{C} \rightarrow \mathfrak{X}
    $$
    from a twisted curve, satisfying the usual stability condition on the induced map of coarse spaces. The evaluation maps of twisted stable maps land in \(\overline{\mathcal{I}}_{\mu}\mathfrak{X}\).
\end{definition}

\begin{definition}[The moduli stack of twisted stable maps]\label{def:twisted-map-moduli}
    Let $\mathfrak{X}$ be an orbifold with coarse moduli space $X$. Let $g,n$ be non-negative integers and $\beta$ be a curve class on $X$. The \emph{moduli stack of twisted stable maps to $\mathfrak{X}$ of genus $g$ curves of class $\beta \in H_{2}(X,\mathbb{Z})$ with $n$ marked points} (we denote it by $\mathcal{K}_{g,n}(\mathfrak{X},\beta)$) is the Deligne-Mumford stack of twisted stable maps $f: \mathfrak{C} \rightarrow \mathfrak{X}$ of curves of genus $g$ with $n$ marked points such that the underlying map satisfies $f_{*}C = \beta$. The notation $H_{2}(-,\mathbb{Z})$ denotes singular homology of the underlying complex analytic space. Curve classes for orbifold targets are classes on their coarse moduli spaces. We retain the full integral homology class, without quotienting by torsion or numerical equivalence.
\end{definition}

\begin{definition}[Gromov-Witten invariants with gravitational descendants, cf. \protect{\cite[\S 8.3]{AGV08}}, \protect{\cite[\S 2.1]{Iritani20}}]\label{definition: Gromov-Witten for orbifold}

    Notation as in \cref{def:twisted-map-moduli}. On $\mathcal{K}_{g,n}(\mathfrak{X},\beta)$, one defines $n$ tautological line bundles $\mathcal{L}_{i}$ whose fibre at a point is the cotangent space to the corresponding \emph{coarse} curve at the $i$-th marked point. The \emph{cotangent line class} is
    $$
    \psi_{i}=c_{1}(\mathcal{L}_{i}) \in H^{2}\!\left(\mathcal{K}_{g,n}(\mathfrak{X},\beta)\right).
    $$
    
    Let $\overline{\mathcal{I}}\mathfrak{X}$ be the rigidified cyclotomic inertia stack of $\mathfrak{X}$. We consider the evaluation maps 
    \begin{align*}
        ev_{i}: \mathcal{K}_{g,n}(\mathfrak{X},\beta) &\rightarrow \overline{\mathcal{I}}\mathfrak{X},\\
        ev_{i}(C,\mu;p_{1}, \cdots , p_{n},f) &= (f(p_{i}),\mu).
    \end{align*}
    Consider $\gamma_{1}, \cdots, \gamma_{n} \in H^{*}(\overline{\mathcal{I}}\mathfrak{X})$. Let $a_{1},\cdots, a_{n}$ be non-negative integers. Then the \emph{Gromov-Witten invariant with descendants} is the number given by
    $$
    \langle \tau_{a_{1}}\gamma_{1}, \cdots , \tau_{a_{n}}\gamma_{n} \rangle_{\beta} = \int_{[\mathcal{K}_{g,n}(\mathfrak{X},\beta)]^{virt}} \psi_{1}^{a_{1}}ev_{1}^{*}(\gamma_{1}) \cdots \psi_{n}^{a_{n}}ev_{n}^{*}(\gamma_{n}).
    $$
\end{definition}

\begin{convention}[Notation for the smooth case]
   When $X$ is a smooth variety, the moduli stack of stable maps will be denoted by $\overline{M}_{g,n}(X,\beta)$.
\end{convention}

\begin{statement}[Set-up]\label{setup:I-series}
        During this subsection, we set up the following data:
    \begin{enumerate}
        \item Let $\mathfrak{X}$ be the canonical covering stack of a Fano threefold $X$ with terminal quotient singularities. Let $n$ be the dimension of $X$ and $r$ be the Picard number of $X$.
        \item Let $D$ be a $\mathbb{C}$-divisor on $X$.
        \item Let $K\subseteq H_{2}(X,\mathbb{Z})$ be the set of curve classes for which the evaluation morphism $\mathrm{ev}_{1}:\mathcal{K}_{0,1}(\mathfrak{X},\beta)\rightarrow X$ is surjective.
        \item Let 
        $$
        \hat{G}_{X}(t) = \widetilde{I}^{X}_{0}(t) = 1 + \sum\limits_{\substack{\beta \in K \\ -K_{X}\cdot\beta \geq 2}}(-K_{X} \cdot \beta)!\langle \tau_{-K_{X}\cdot\beta -2}\mathbf{1} \rangle_{\beta} \cdot t^{-K_{X} \cdot \beta}
        $$
        be \emph{the constant term of regularized $I$-series} or \emph{the regularized quantum period} of $X$. Here the class $\mathbf{1} \in H^{2 \mathrm{dim}X}(\overline{\mathcal{I}}\mathfrak{X},\mathbb{Q})$ is the Poincaré dual of the point class of the untwisted sector of $\overline{\mathcal{I}}\mathfrak{X}$.
        \item Let 
        \begin{align*}
        \hat{G}_{X,D}(t) = \widetilde{I}_{0}^{X,D}(t) &= 1 + \sum\limits_{\substack{\beta \in K \\ -K_{X}\cdot\beta \geq 2}}(-K_{X} \cdot \beta)!\langle \tau_{-K_{X}\cdot\beta -2}\mathbf{1} \rangle_{\beta} \cdot e^{-D \cdot \beta}t^{-K_{X} \cdot \beta}\\
        &= 1 + a_{1}t + a_{2}t^2 + \cdots
        \end{align*} 
        be \emph{the restriction of the constant term of regularized $I$-series of $X$ to the
        anti-canonical direction corresponding to $D$}.
    \end{enumerate}
\end{statement}

\begin{remark}
    The above definition of the quantum period is a special case of the definition given in \cite{OnetoPetracci2018}.
\end{remark}
    
    Now we recall the orbifold degeneration formula. We start with the definition of transversal intersections for orbifolds.

\begin{definition}[Transversal intersection of smooth Deligne-Mumford stacks, cf. \protect{\cite[Appendix A]{AF16}}]

\hspace{2em}

\begin{itemize}
    \item Let $X$ be a complex algebraic stack. We say that $X$ has \emph{nodal (codimension 1) singularities} if it is locally isomorphic in the f.p.p.f. topology to $\mathrm{Spec}\frac{\mathbb{C}[x,y,z_{1},\cdots,z_{n}]}{(xy)}$. In particular, the singular locus $D$ of $X$ is a smooth substack. Let $\nu: \widetilde{X} \rightarrow X$ be the normalization of $X$ and $\widetilde{D} = \nu^{-1}D$.
    \item A morphism $f: C \rightarrow X$ between nodal algebraic stacks is called \emph{transversal to the singular locus} if
    \begin{enumerate}
        \item the induced morphism $\widetilde{C} \rightarrow \widetilde{X}$ defines a morphism of locally smooth pairs which is transversal to the boundary divisor;
        \item for every point $p \in f^{-1}D$ its two inverse images in $\widetilde{C}$ map to different points of $\widetilde{D}$ via $f$.
    \end{enumerate}
\end{itemize}
\end{definition}

    We also need the relative orbifold Gromov-Witten theory for the orbifold degeneration formula. First, we set up the data we need to define the relative orbifold Gromov-Witten theory.

\begin{statement}[Data for a pair]\label{data for a pair}
    Fix a smooth orbifold pair $(X,D)$. We also fix
    \begin{enumerate}
        \item an integral curve class $\beta\in H_{2}(X,\mathbb{Z})$;
        \item an integer $g \geq 0$;
        \item disjoint finite ordered sets $N,M$ which we may take to be $\{1,2,\cdots,n \}$ and $\{n+1,n+2,\cdots,n + |M|\}$;
        \item tuples $\mathbf{e} = (e_{i})_{i \in N}$ and $\mathbf{f} = (f_{j})_{j \in M}$ of positive integers such that $\mathcal{I}_{e_{i}}(X) \neq \emptyset$ for all $i \in N$ and $\mathcal{I}_{f_{j}}(D) \neq \emptyset$ for all $j \in M$.
        \item a tuple $\mathbf{c} = (c_{j})_{j \in M}$ of positive integers such that $\sum\limits_{j \in M}\frac{c_{j}}{f_{j}} = D \cdot \beta$.
    \end{enumerate}
    We denote by $\boldsymbol\mu = (\mu_{j})_{j \in M}$ the tuple formed by $\mu_{j} = \frac{c_{j}}{f_{j}}$.
\end{statement}

\begin{definition}[Relative Gromov-Witten invariants with gravitational descendants, cf. {\cite[Definition 4.5]{AF16}}]
    Notation as in \cref{data for a pair}. Put
    $$
    \Gamma=(g,\beta,N,M,\mathbf{e},\mathbf{f},\mathbf{c}).
    $$    
    Fix a twisting choice $\mathfrak{r}$, and write $\mathcal{K}_{\Gamma}^{\mathfrak{r}}(X,D)$ for the relative moduli stack with these data. Let $\psi_{i}$ and $\mathrm{ev}_{i}$ denote its coarse cotangent-line classes and evaluation maps. For integers $m_{i}\geq 0$, $i\in N$, and classes
    $$
    \gamma_{i}\in H^{*}(\overline{\mathcal{I}}(X),\mathbb{Q}) \quad(i\in N), \qquad \gamma_{j}\in H^{*}(\overline{\mathcal{I}}(D),\mathbb{Q}) \quad(j\in M),
    $$
    define the relative descendant invariant by
    $$
    \langle \prod\limits_{i \in N} \tau_{m_{i}}(\gamma_{i}) \mid \prod\limits_{j \in M} \gamma_{j} \rangle_{\Gamma}^{(X,D)} := \mathrm{deg}((\prod\limits_{i \in N}\psi_{i}^{m_{i}}\cdot \mathrm{ev}_{i}^{*}\gamma_{i}) \cdot (\prod\limits_{j \in M}\mathrm{ev}_{j}^{*}\gamma_{j}) \cap [\mathcal{K}_{\Gamma}^{\mathfrak{r}}(X,D)]^{vir}).
    $$
\end{definition}

    Now we are ready to state the orbifold degeneration formula.

\begin{convention}[Data for a degeneration]\label{data for a degeneration}
    Fix a flat morphism $\pi: W \rightarrow B$ such that $B$ is a smooth curve, $W$ is a smooth Deligne-Mumford stack, and $0 \in B$ is the unique critical value of $\pi$. We assume that $W_{0} = X_{1} \cup_{D} X_{2}$ is the union of two smooth closed substacks $X_{1}$ and $X_{2}$ intersecting transversally along $D$, a smooth divisor on both $X_{1}$ and $X_{2}$. This implies that $W_{0}$ is nodal and first-order smoothable along the singular locus $D$. The rigidified inertia stack of $W_{0}$ is denoted by $\overline{\mathcal{I}}(W_{0})$. We also fix
    \begin{enumerate}
        \item a curve class $\beta$ in the fibre of $\pi$;
        \item an integer $g \geq 0$;
        \item a finite ordered set $N$, possibly empty, which we may take to be $\{1,2,\cdots,n \}$;
        \item a tuple $\mathbf{e} = (e_{i})_{i \in N}$ of positive integers such that $\mathcal{I}_{e_{i}}(W_{b}) \neq \emptyset$ for all $b \in B, i \in N$. 
    \end{enumerate}
\end{convention}

\begin{theorem}[The orbifold degeneration formula, {\cite[\S 5]{AF16}}]\label{theorem: orbifold degeneration formula}
    Notation as in \cref{data for a degeneration}. Assume that $W_{0}$ is proper with projective coarse moduli space. Let
    $$
    \gamma_1,\ldots,\gamma_n \in H^{*}(\overline{\mathcal{I}}(W_{0}),\mathbb{Q})
    $$
    have homogeneous parity, and let $a_1,\ldots,a_n$ be nonnegative integers. Then
    \begin{align*}
        \langle \prod\limits_{i=1}^{n}\tau_{a_{i}}(\gamma_{i}) \rangle^{W_{0}}_{g,\beta} = \sum\limits_{\eta \in \Omega} \frac{\prod\limits_{j \in M}d_{j}}{|M|!}\sum\limits_{\substack{\delta_{j} \in F \\ j \in M}}(-1)^{\varepsilon} \langle \prod\limits_{i \in N_{1}}\tau_{a_{i}}(\gamma_{i}) \mid \prod\limits_{j \in M} \delta_{j} \rangle^{(X_{1},D)}_{\Xi_{1}} \cdot \langle \prod\limits_{i \in N_{2}}\tau_{a_{i}}(\gamma_{i}) \mid \prod\limits_{j \in M} \widetilde{\delta}^{\vee}_{j} \rangle^{(X_{2},D)}_{\Xi_{2}}
    \end{align*}
    where
    \begin{enumerate}
        \item $F$ is a homogeneous basis of $H^{*}(\overline{\mathcal{I}}(D),\mathbb{Q})$.
        \item $\widetilde{\delta}^{\vee}_{j}$ is the dual of $\delta_{j}$ with respect to the Chen-Ruan pairing.
        \item $\Omega$ is the set of splittings of the data $g,n,\beta$. An element in $\Omega$ includes the data below:
        \item $N_{1},N_{2}$ is a decomposition of $\{1,2,\cdots,n \}$ in two subsets.
        \item $\Xi_{1},\Xi_{2}$ is a possibly disconnected splitting of the data $\beta,g$ in two modular graphs having roots labeled by $M:=\{n+1,n+2,\cdots,n+|M|\}$.
        \item $d_{i}$ are assigned intersection multiplicities for $i \in M$.
        \item $(-1)^{\varepsilon}$ is the sign determined formally by the equality
        $$
        \prod\limits_{i=1}^{n} \gamma_{i} \prod\limits_{j \in M}\delta_{j}\widetilde{\delta}^{\vee}_{j} = (-1)^{\varepsilon} \prod\limits_{i \in N_{1}} \gamma_{i} \prod\limits_{j \in M}\delta_{j} \prod\limits_{i \in N_{2}}\gamma_{i} \prod\limits_{j \in M} \widetilde{\delta}^{\vee}_{j}.
        $$
    \end{enumerate}
\end{theorem}

    We recall what it means for a smooth center to have normal crossings with a boundary.

\begin{definition}[Normal crossings with a boundary]
    Let $\mathfrak{Y}$ be a smooth algebraic stack, let $\mathfrak{E}=\sum\limits_{i\in I}\mathfrak{E}_{i}$ be a labeled reduced simple-normal-crossings divisor, and let $\mathfrak{C}\subseteq\mathfrak{Y}$ be a closed substack. The following conditions are checked after pullback to a smooth atlas $U\rightarrow\mathfrak{Y}$. Write $C_U$ and $E_U$ for the corresponding pullbacks.
    
    We say that $\mathfrak{C}$ \emph{has normal crossings with} $\mathfrak{E}$, or \emph{meets $\mathfrak{E}$ normally}, if \'{e}tale locally at every point $P\in C_U$ there are regular parameters $x_{1},\ldots,x_{n}$ such that
    $$
    I_{C_U,P}=(x_{1},\ldots,x_{m}),
    $$
    and every irreducible component of $E_U$ through $P$ is defined by one of these parameters.
    
    If, in addition, the parameters defining the boundary components through $P$ can all be chosen among $x_{m+1},\ldots,x_{n}$, then $\mathfrak{C}$ meets every stratum of $\mathfrak{E}$ transversely. In this case we say that $\mathfrak{C}$ \emph{is transverse to the boundary}. For a Deligne-Mumford stack, the atlas may be chosen \'{e}tale.
\end{definition}

\subsection{The log smooth degeneration formula}\label{subsec:scheme-log-degeneration}

    In this subsection, we recall the logarithmic degeneration formula of \cite{ACGS20} and the component-wise cutting-and-gluing theorem of \cite{Ran22}.

\begin{definition}[Semistable morphisms]\label{def:semistable-morphism}
    Let $0 \in B$ be a pointed smooth curve. We work after shrinking $B$ around $0$ when necessary. A morphism $\pi:Y\rightarrow B \ni 0$ is \emph{semistable} if $Y$ is smooth, $\pi$ is smooth over $B\setminus\{0\}$, and the central fibre $Y_{0}$ is a reduced simple-normal-crossings divisor on $Y$.
\end{definition}

\begin{definition}[Log smooth degenerations]\label{def:scheme-log-degeneration}
    Let $0 \in B$ be a pointed smooth curve, and let $\pi: Y\rightarrow B \ni 0$ be a projective semistable morphism. Endow $B$ with the divisorial logarithmic structure associated with the divisor $0$ and endow $Y$ with the divisorial logarithmic structure associated with the divisor $Y_{0}$. Denote the resulting logarithmic schemes by $B^{\dagger}$ and $Y^{\dagger}$. Then the morphism $\pi$ induces a logarithmic morphism
    $$
    \pi^{\dagger} : Y^{\dagger} \rightarrow B^{\dagger}.
    $$
    We say that $Y^{\dagger} \rightarrow B^{\dagger}$ is a \emph{log smooth degeneration}.
\end{definition}

\begin{remark}
    The results of \cite{ACGS20} apply more generally to projective logarithmically smooth morphisms whose special fibre is simple in the sense of \cite[Definition 2.1]{ACGS20}. The semistable situation above is sufficient for all applications in this article.
\end{remark}

\begin{definition}[Stable logarithmic maps]\label{def:stable-logarithmic-maps}
    Let $Y^{\dagger}\rightarrow B^{\dagger}$ be a projective log smooth degeneration, and let $S^{\dagger}\rightarrow B^{\dagger}$ be a morphism from a fine and saturated logarithmic scheme. A \emph{logarithmic curve over $S^{\dagger}$} is a logarithmically smooth and integral morphism
    $$
    C^{\dagger}\rightarrow S^{\dagger}
    $$
    whose underlying morphism is a family of nodal marked curves.
    
    A \emph{stable logarithmic map to $Y^{\dagger}/B^{\dagger}$ over $S^{\dagger}$} is a logarithmic morphism
    $$
    f^{\dagger}: C^{\dagger} \rightarrow Y^{\dagger}
    $$
    such that the underlying morphism
    $$
    f: C \rightarrow Y
    $$
    is a stable map over $B$. A stable logarithmic map is called \emph{basic} if its logarithmic structure on the base is minimal among the logarithmic structures supporting the given underlying map and its contact data.

    For a stable logarithmic map, we have the associated data $\Gamma=(g,n,\beta,\mathbf{u})$ where:
    \begin{itemize}
        \item $g$ is the genus;
        \item $n$ is the number of marked points;
        \item $\beta\in H_{2}(Y,\mathbb{Z})$ is the pushforward of the fibre curve class to the total space;
        \item $\mathbf{u}=(u_{1},\ldots,u_{n})$ is the tuple of prescribed contact vectors at the marked points.
    \end{itemize}

    Denote by $\mathscr{M}_{\Gamma}^{\log}(Y/B)$ the moduli stack of basic stable logarithmic maps with these data. For $b\in B$, let $\iota_b:Y_b\hookrightarrow Y$ be the inclusion. The fibre mapping stack $\mathscr{M}_{\Gamma}^{\log}(Y_b/b)$ includes the integral classes $\beta_b\in H_{2}(Y_b,\mathbb{Z})$ satisfying
    $$
    (\iota_b)_{*}\beta_b=\beta \quad\text{in }H_{2}(Y,\mathbb{Z}).
    $$
    We use the integral homology option for degree data in \cite[\S~2.3.8]{ACGS20}.
\end{definition}

\begin{proposition}[The logarithmic deformation invariance, cf. {\cite[Theorem 1.1]{ACGS20}}]\label{prop: log deformation invariance}
    Notation as in \cref{def:stable-logarithmic-maps}, with $\Gamma$ fixed. Then for the inclusion $j_b:b\hookrightarrow B$ of a closed point, we have
    $$
    j_{b}^{!} \left[
    \mathscr{M}_{\Gamma}^{\log}(Y/B) \right]^{\mathrm{vir}} = \left[
    \mathscr{M}_{\Gamma}^{\log}(Y_{b}/b) \right]^{\mathrm{vir}}.
    $$
    In particular, for descendant insertions obtained by restricting cohomology classes on $Y$, the corresponding generic-fibre and central-fibre virtual integrals agree.
\end{proposition}

\begin{definition}[Tropical types]\label{definition: tropical types}
    Notation as in \cref{prop: log deformation invariance}. Write $\Sigma(Y_{0})$ for the cone complex of the logarithmic central fibre. In a decorated type below, $i_{v}$ denotes the inclusion into $Y_{0}$ of the closure of the target stratum assigned to $v$, and $A(v)$ is an integral curve class on that closure. A \emph{tropical type} $\tau = (G,g,\sigma,\mathbf{u})$ consists of:
    \begin{itemize}
        \item a connected graph $G$ with legs;
        \item a genus $g(v) \in \mathbb{Z}_{\geq 0}$ for every vertex $v$;
        \item a cone $\sigma(x) \subset \Sigma(Y_{0})$ for every vertex, edge and leg $x$ of $G$;
        \item an integral contact vector $u_{h}$ for every oriented edge or leg $h$, subject to the incidence and balancing conditions, and to
        $$
        u_{\overline{h}} = -u_{h}
        $$
        for the two orientations of every bounded edge.
    \end{itemize}

    A \emph{decorated tropical type}
    $$
    \widetilde{\tau} = (\tau, A)
    $$
    is a tropical type together with an effective integral curve class $A(v)$ on the stratum closure assigned to each vertex. Let $\iota_{0}:Y_{0}\hookrightarrow Y$ be the inclusion. The degree condition is
    $$
    \beta = \sum\limits_{v\in V(G)}(\iota_{0}\circ i_v)_{*}A(v) \quad\text{in }H_{2}(Y,\mathbb{Z}).
    $$

    A tropical type is called \emph{rigid} if its cone of tropical realizations is one-dimensional and maps nontrivially to the cone of the logarithmic base.
\end{definition}

\begin{proposition}[The logarithmic degeneration formula, cf. {\cite[Theorem 1.2]{ACGS20}}]\label{prop: scheme log degeneration formula}
    Notation as in \cref{prop: log deformation invariance}. For a decorated rigid tropical type $\widetilde{\tau}$, let $\mathscr{M}^{\log}(Y_{0},\widetilde{\tau})$ be the moduli stack of maps marked by $\widetilde{\tau}$, and let
    $$
    j_{\widetilde{\tau}}: \mathscr{M}^{\log}(Y_{0},\widetilde{\tau}) \rightarrow \mathscr{M}_{\Gamma}^{\log}(Y_{0}/0)
    $$ 
    forget that marking. For the underlying rigid tropical type $\tau$, let $Q_\tau$ be its basic monoid. The multiplicity $m_\tau$ is the image of the primitive generator of the dual ray $Q_\tau^\vee$ in the integral lattice of the base ray (cf. \cite[Definition~3.6 and Corollary~3.8]{ACGS20}). Then
    $$
    \left[
    \mathscr{M}_{\Gamma}^{\log}(Y_{0}/0)
    \right]^{\mathrm{vir}} 
    = \sum\limits_{\widetilde{\tau}} \frac{m_{\tau}}{|\mathrm{Aut}(\tau)|}(j_{\widetilde\tau})_{*} \left[
    \mathscr{M}^{\log}(Y_{0},\widetilde\tau)
    \right]^{\mathrm{vir}}
    $$
    where the sum runs over all decorations of a fixed representative of every isomorphism class of rigid tropical types of total discrete data $\Gamma$. The group $\mathrm{Aut}(\tau)$ consists of automorphisms of the graph preserving its legs, genera, target cones, and contact vectors.
\end{proposition}

    The expanded theory below describes the individual summands in \cref{prop: scheme log degeneration formula} through relative mapping problems for the components of the special fibre.

\begin{definition}[Stable logarithmic maps to expansions]\label{def:expanded-logarithmic-maps}
    Notation as in \cref{prop: scheme log degeneration formula}. We write $\widetilde{\tau}$ for the decorated rigid tropical type denoted by $\rho$ in \cite{Ran22}. Let $\lambda$ be a polyhedral choice in the inverse system of refinements of the universal tropical map, and let
    $$
    Y^{\lambda} \rightarrow Y
    $$
    be the associated logarithmic expansion. A \emph{stable logarithmic map to a $\lambda$-expansion of $Y$} is a diagram
    $$
    C^{\lambda} \rightarrow Y^{\lambda} \rightarrow Y
    $$
    such that $C^{\lambda}\rightarrow Y^{\lambda}$ is a stable logarithmic map and its tropicalization is combinatorially transverse: every vertex maps to a vertex and every non-contracted edge maps onto an edge. We denote by $K_{\Gamma}^{\lambda}(Y)$ the moduli stack of stable logarithmic maps to $\lambda$-expansions of $Y$ with discrete data $\Gamma$. If $Y=Y_{0}$, we denote by $K_{\widetilde{\tau}}^{\lambda}(Y_{0})$ the moduli stack of such maps together with a marking of the unexpanded stabilization of
    $$
    C^{\lambda}\rightarrow Y^{\lambda}\rightarrow Y_{0}
    $$
    by the decorated rigid tropical type $\widetilde{\tau}$, in the sense of \cite[Definition 2.31]{ACGS20}. At each geometric point, the decorated tropical type of the unexpanded stabilization is equipped with a contraction to $\widetilde{\tau}$, compatibly with generalization (cf. \cite[Remark 2.32]{ACGS20}).
\end{definition}

\begin{proposition}[The expanded cutting-and-gluing, cf. {\cite[Theorem A and B, Lemma 7.4.3-7.4.4]{Ran22}}]\label{prop: expanded cut and glue}
    Notation as in \cref{def:expanded-logarithmic-maps}. Fix $\Gamma$, a decorated rigid marking type $\widetilde{\tau}$, and a sufficiently fine polyhedral choice $\lambda$. Let $G$ be the graph of the chosen transverse refinement. For $v\in V(G)$, let $Y_v$ be the corresponding component of the normalization of the modified central fibre, and let $K_{\widetilde{\tau}}^{\lambda}(Y_v)$ be its full expanded relative mapping stack with the prescribed markings, curve class, and contact orders. Further compatible logarithmic refinements are permitted.

    Then the following statements hold:
    \begin{enumerate}
        \item \emph{Comparison with stable logarithmic maps.} The stack $K_{\Gamma}^{\lambda}(Y_{0})$ is a proper Deligne–Mumford stack equipped with a virtual fundamental class. There is a canonical logarithmic modification
        $$
        p^{\lambda}: K_{\Gamma}^{\lambda}(Y_{0}) \rightarrow \mathscr{M}_{\Gamma}^{\log}(Y_{0}/0)
        $$
        such that
        $$
        p^{\lambda}_{*}[K_{\Gamma}^{\lambda}(Y_{0})]^{\mathrm{vir}} = 
        \left[
        \mathscr{M}_{\Gamma}^{\log}(Y_{0}/0)
        \right]^{\mathrm{vir}}.
        $$
        In particular, the logarithmic Gromov–Witten invariants of $Y_{0}$ may be computed using stable logarithmic maps to expansions.
        \item \emph{The modified product of vertex theories.} Cutting the bounded edges of $G$ gives a morphism
        $$
        K_{\widetilde{\tau}}^{\lambda}(Y_{0}) \rightarrow \prod_{v\in V(G)} K_{\widetilde{\tau}}^{\lambda}(Y_{v}).
        $$
        There is a logarithmic modification
        $$
        \bigtimes\limits_{v\in V(G)}^{\lambda}
        K_{\widetilde{\tau}}(Y_{v})
        \longrightarrow
        \prod_{v\in V(G)}
        K_{\widetilde{\tau}}^{\lambda}(Y_{v})
        $$
        called the modified product, and, after replacing $\lambda$ by a sufficiently fine refinement, the cutting morphism lifts to
        $$
        \kappa:
        K_{\widetilde{\tau}}^{\lambda}(Y_{0})
        \longrightarrow
        \bigtimes\limits_{v\in V(G)}^{\lambda}
        K_{\widetilde{\tau}}(Y_{v}).
        $$
        The modified product carries a natural virtual fundamental class
        $$
        \left[
        \bigtimes\limits_{v\in V(G)}^{\lambda}
        K_{\widetilde{\tau}}(Y_{v})
        \right]^{\mathrm{vir}}.
        $$
        \item \emph{The relative diagonal and the gluing space.} For every bounded edge $e\in E(G)$, let $D_{e}(Y^{\lambda})$ be the family of open divisors supporting the evaluation morphisms at the two marked points obtained by cutting $e$. The two flags of $e$ determine an evaluation morphism
        $$
        \mathrm{ev}_{e}:
        \bigtimes\limits_{v\in V(G)}^{\lambda}
        K_{\widetilde{\tau}}(Y_{v})
        \longrightarrow
        D_{e}^{\{2\}}(Y^{\lambda}),
        $$
        where $D_{e}^{\{2\}}(Y^{\lambda})$ is the relative fibre square $D_{e}(Y^{\lambda}) \times D_{e}(Y^{\lambda})$. Let $D_{G}$ be the fibre product of the divisor families $D_{e}(Y^{\lambda})$ over all bounded edges $e$ of $G$, and denote by $D_{G}^{2}$ its relative fibre square. The edge evaluation morphisms combine to give
        $$
        \mathrm{ev}:
        \bigtimes\limits_{v\in V(G)}^{\lambda}
        K_{\widetilde{\tau}}(Y_{v}) \rightarrow D_{G}^{2}.
        $$
        Let $\Delta:D_{G}\longrightarrow D_{G}^{2}$ be the relative diagonal. The gluing space $G^{\lambda}(Y_{0})$ is defined by the Cartesian diagram
        $$
        \begin{tikzcd}
        G^{\lambda}(Y_{0}) \arrow[r,"\mathrm{split}"] \arrow[d] & \bigtimes\limits_{v\in V(G)}^{\lambda} K_{\widetilde{\tau}}(Y_{v}) \arrow[d,"\mathrm{ev}"] \\
        D_{G} \arrow[r,"\Delta"] & D_{G}^{2}
        \end{tikzcd}
        $$
        Equivalently,
        $$
        G^{\lambda}(Y_{0}) = \left(
        \bigtimes\limits_{v\in V(G)}^{\lambda}
        K_{\widetilde{\tau}}(Y_{v}) \right) \times_{D_{G}^{2}} D_{G}.
        $$
        Its virtual fundamental class is
        $$
        [G^{\lambda}(Y_{0})]^{\mathrm{vir}} = \Delta^{!} \left[
        \bigtimes\limits_{v\in V(G)}^{\lambda}
        K_{\widetilde{\tau}}(Y_{v})
        \right]^{\mathrm{vir}}.
        $$
        
        The space $G^{\lambda}(Y_{0})$ parametrizes collections of maps to expansions of the components $Y_{v}$ whose underlying maps glue along the pairs of marked points corresponding to the bounded edges of $G$.
        \item \emph{Logarithmic structures at the gluing nodes.}
        After the reductions and compatible refinements of \cite[\S 7.4]{Ran22}, there is a finite \'etale morphism
        $$
        \mu: K_{\widetilde{\tau}}^{\lambda}(Y_{0}) \rightarrow G^{\lambda}(Y_{0}).
        $$
        Its obstruction theory is pulled back from the gluing space, and
        $$
        [K_{\widetilde{\tau}}^{\lambda}(Y_{0})]^{\mathrm{vir}} = \mu^{*}[G^{\lambda}(Y_{0})]^{\mathrm{vir}}.
        $$
        On a component on which $\mu$ has constant degree,
        $$
        \mu_{*}[K_{\widetilde{\tau}}^{\lambda}(Y_{0})]^{\mathrm{vir}} = \deg(\mu)[G^{\lambda}(Y_{0})]^{\mathrm{vir}}.
        $$
        The degree and all tropical multiplicities are taken with the same logarithmic-base conventions.
    \end{enumerate}
\end{proposition}

\section{Comparison lemmas for degenerations}\label{sec:scheme-comparisons}

    In this section, we establish some comparison lemmas of Gromov-Witten invariants that will be used in \cref{sec:stationary-terminal-invariants} and \cref{sec:ordinary-terminal-case}. We use logarithmic deformation invariance, tropical decomposition, and expanded cutting and gluing as recalled in \cref{prop: log deformation invariance}, \cref{prop: scheme log degeneration formula}, and \cref{prop: expanded cut and glue}.

\subsection{Virtual dimension and vanishing before gluing}\label{subsec:scheme-virtual-dimension}

\begin{proposition}[Virtual dimension of a logarithmic factor]\label[proposition]{prop:relative-virtual-dimension}
    Let $Z$ be a smooth projective threefold and $B=\sum\limits_{i=1}^{r}B_{i}$ be a reduced simple normal crossings divisor. Consider connected stable logarithmic maps to $(Z,B)$ of genus $g$ and class $\alpha$, with $n$ original markings of contact vector zero and $l$ relative markings. Let $\mu_{h,i}\in\mathbb{Z}_{\geq 0}$ be the contact order with $B_{i}$ at the $h$-th relative marking. Denote the moduli stack for these fixed data by $\mathcal{M}$. Then
    $$
    B_{i}\cdot\alpha = \sum\limits_{h=1}^{l}\mu_{h,i}
    $$
    and
    $$
    \begin{aligned}
    \mathrm{vdim}\mathcal{M}
    &=
    (-K_{Z}-B)\cdot\alpha+n+l\\
    &=
    (-K_{Z})\cdot\alpha+n+l
    -\sum\limits_{h=1}^{l}\sum\limits_{i=1}^{r}\mu_{h,i}.
    \end{aligned}
    $$
    The same virtual dimension is used on the compatible expanded relative factors.
\end{proposition}

\begin{proof}
    The logarithmic tangent bundle $T_{Z}(-\log B)$ has first Chern class
    $$
    c_{1}\!\left(T_{Z}(-\log B)\right)=-K_{Z}-B.
    $$
    By Riemann-Roch and the contribution from the marked source curve, we have
    $$
    \begin{aligned}
    \mathrm{vdim}\mathcal{M}
    &=
    (1-g)(3-3)
    +(-K_{Z}-B)\cdot\alpha+n+l\\
    &=
    (-K_{Z}-B)\cdot\alpha+n+l.
    \end{aligned}
    $$

    For the contact identity, normalize the source curve and sum the degrees of the pullback of $\mathcal{O}_{Z}(B_{i})$ over its components. The signed contributions from the two branches of each internal node cancel. The remaining contributions are the marking contacts. Hence
    $$
    B_{i}\cdot\alpha=\sum\limits_h\mu_{h,i}.
    $$
    The expanded obstruction theories and their compatible refinements preserve the virtual rank (cf. \cite[\S~4.5 and Proposition~4.6.1]{Ran22}).
\end{proof}

\begin{lemma}[Vanishing of a relative factor]\label{lem:scheme-factor-vanishing}
    Notation as in \cref{prop: expanded cut and glue}. Fix one decorated contribution and the individual compatible refinements. For $v\in V(G)$, write $\mathcal{K}_v$ for its relative mapping stack. Let $v_{0}\in V(G)$, and suppose that an insertion $\eta_{v_{0}}$ on $\mathcal{K}_{v_{0}}$ satisfies
    $$
    \eta_{v_{0}}\cap [\mathcal{K}_{v_{0}}]^{\mathrm{vir}}=0.
    $$
    Then the entire glued contribution containing the pullback of this insertion is zero.
\end{lemma}

\begin{proof}
    We use the virtual constructions of \cite[\S\S~4.4-4.6 and~7.3-7.4]{Ran22}. Write
    $$
    P=\prod_{v}\mathcal{K}_{v}.
    $$
    After regular refinements of the combinatorial bases, the factors have strict morphisms to smooth combinatorial mapping stacks $S_{v}$. Put $S=\prod_{v}S_{v}$.
    
    The subdivision synchronizing the expansions induces a morphism
    $$
    b:S^{\lambda}\rightarrow S
    $$
    between smooth combinatorial stacks. The modified product is the ordinary Cartesian pullback
    $$
    P^{\lambda}=P\times_{S} S^{\lambda}.
    $$
    Strictness identifies this with the required logarithmic pullback. The relative obstruction theory is pulled back from that of $P/S$.
    
    The morphism $b$ is lci of virtual relative dimension zero. Compatibility of the obstruction theories gives
    $$
    [P^{\lambda}]^{\mathrm{vir}}
    =
    b_P^{!}
    \left(
    \mathop{\boxtimes}_{v}
    [\mathcal{K}_{v}]^{\mathrm{vir}}
    \right).
    $$
    Let
    $$
    \eta=\mathrm{pr}_{v_{0}}^{*}\eta_{v_{0}}.
    $$
    The external-product formula gives
    $$
    \eta\cap
    \left(
    \mathop{\boxtimes}_{v}
    [\mathcal{K}_{v}]^{\mathrm{vir}}
    \right)=0.
    $$
    Compatibility of refined pullback with insertions shows that its pullback annihilates $[P^{\lambda}]^{\mathrm{vir}}$.
    
    The relative diagonal pullback, the finite \'etale enhancement pullback, and the final proper pushforward preserve this zero class. Thus the complete contribution vanishes.
    
    For cohomological insertions, the same argument is applied to the cohomological realization of the virtual classes and refined operations.
\end{proof}

\subsection{The discrepancy calculation}\label{subsec:scheme-discrepancy}

\begin{lemma}[Discrepancies and contact orders]\label[lemma]{lem:ordinary-discrepancy-contact}
    Let $X$ be a terminal threefold, $D$ be an effective $\mathbb{Q}$-Cartier divisor on $X$, and 
    $$
    p: Z \rightarrow X
    $$ 
    be a log resolution of $(X,D)$. Let $B=\sum\limits_{i=1}^{r}B_{i}$ be the sum of the exceptional divisors of $p$ and the strict transformations of the prime components of $D$. Write
    $$
    K_{Z}=p^{*}K_{X}+\sum\limits_{i=1}^{r} a_{i}B_{i}, \qquad p^{*}D=\sum\limits_{i=1}^{r} m_{i}B_{i}.
    $$
    Then, for $c>0$ and a nonzero vector $(\mu_{i})_{i}\in\mathbb{Z}_{\geq 0}^{\,r}$, we have
    $$
    \sum\limits_{i=1}^{r}(a_{i}+1+cm_{i})\mu_{i}>1.
    $$
    Moreover, if $D=0$, then the same assertion holds for $c=0$.
\end{lemma}

\begin{proof}
    If $B_{i}$ is exceptional over $X$, then $a_{i}>0$. If $B_{i}$ is not exceptional, then $a_{i}=0$ and $m_{i}>0$. In particular, we have
    $$
    (a_{i}+1+cm_{i})\mu_{i} \geq 0
    $$
    for all $1 \leq i \leq r$. Choose $i$ such that $\mu_{i}>0$. Since $\mu_{i}$ is an integer, $\mu_{i}\geq 1$. If $B_{i}$ is exceptional, then
    $$
    (a_{i}+1)\mu_{i} > 1.
    $$
    If $B_{i}$ is nonexceptional, then
    $$
    (1+cm_{i})\mu_{i} > 1.
    $$
\end{proof}

\begin{lemma}[Vanishing for a point descendant]\label{lem:ordinary-marked-factor}
    Notation as in \cref{lem:ordinary-discrepancy-contact}. Let $\pi:\mathcal{W}\rightarrow T$ be a projective semistable degeneration over a smooth pointed curve $0\in T$, with distinguished component $Z$ and reduced vertical boundary $B$. Let $\mathcal{L}$ be a relatively nef $\mathbb{Q}$-Cartier divisor on $\mathcal{W}$ such that
    $$
    \mathcal{L}|_{Z} = p^{*}(-K_{X}+cD).
    $$
    Put $L_{X}=-K_{X}+cD$. Assume either that $c>0$, or that $D=0$, $c=0$. Consider connected stable maps with exactly one original marked point, carrying the point insertion. Let this point insertion specialize to $q\in Z\setminus B$. Fix the genus and a total curve-class decoration for connected maps. Let $d \in \mathbb{Z}_{\geq2}$ and let $e \in \mathbb{R}$ such that the total $\mathcal{L}$-degree is $d+ce$. With $\psi$ the coarse cotangent-line class and $\mathrm{ev}$ the evaluation at that marking, consider
    $$
    \psi^{d-2}\mathrm{ev}^{*}\mathbf{1}.
    $$
    If $c>0$ and $e<0$, every contribution vanishes. If $e=0$, every contribution in which the relative factor containing the point has a relative marked point vanishes.
\end{lemma}

\begin{proof}
    Let $\alpha$ be the class of the connected relative factor containing the point, and let $\delta\geq 0$ be the total $\mathcal{L}$-degree of all remaining factors. Then
    $$
    d+ce=p^{*}L_{X}\cdot\alpha+\delta.
    $$
    Let $\mathcal{M}$ be the relative moduli stack of the factor carrying the point insertion. If the $l$ relative markings have contact orders $\mu_{h,i}$, then \cref{prop:relative-virtual-dimension} gives
    $$
    \begin{aligned}
    \mathrm{vdim}\mathcal{M}-(d+1)
    ={}&ce-\delta\\
    &+
    \sum\limits_{h=1}^{l}
    \left(
    1-\sum\limits_{i}(a_{i}+1+cm_{i})\mu_{h,i}
    \right).
    \end{aligned}
    $$
    Each summand in parentheses is strictly negative by \cref{lem:ordinary-discrepancy-contact}.

    The insertion has cohomological degree $2(d+1)$. It annihilates the factor's virtual class when $c>0$ and $e<0$, including the case $l=0$, and when $e=0$ and $l>0$.

    We apply this calculation on the individually refined relative factors. The point remains in the interior, so its descendant pulls back from the marked factor. \cref{lem:scheme-factor-vanishing} then eliminates the entire gluing contribution.
\end{proof}

\begin{lemma}[The contribution with no cut edges]\label[lemma]{lem:ordinary-edge-free}
    Let $\pi:\mathcal{W}\rightarrow T$ be a projective semistable degeneration with distinguished component $Z$ of multiplicity one and reduced vertical boundary
    $$
    B= \left( Z\cap\overline{\mathcal{W}_{0}\setminus Z} \right)_{\mathrm{red}}.
    $$
    Fix discrete data for connected stable maps and their insertions in the degeneration formula of \cref{prop: expanded cut and glue}. Assume that there is at least one original marking, that all nonzero contributions place the original markings on the distinguished component, and that contributions with relative marked points on its factors have vanished.

    Then the remaining contribution is the sum of the logarithmic relative invariants of the distinguished pair with no relative marked points and with the compatible curve classes, each with coefficient one.
\end{lemma}

\begin{proof}
    A cut edge incident to a distinguished factor has nonzero contact. Therefore the surviving factors have no such edges. Connectedness and the presence of an original marking leave the entire source in the distinguished relative theory.

    The contact-degree identity gives zero intersection with every component of its boundary. The surviving marking graph has one vertex and no edges. Its base multiplicity is one because the distinguished component has multiplicity one. There are no diagonal matching conditions, edge-enhancement factors, or graph automorphisms. The result follows from the expanded comparison and its virtual pushforward.
\end{proof}

\section{Stationary Gromov-Witten invariants of terminal Fano threefolds}\label{sec:stationary-terminal-invariants}

    In this section, we construct stationary descendant Gromov-Witten invariants for a general terminal Fano threefold $X$ by good resolutions. The main result in this section is that this definition is well-defined. More precisely, we have
    \begin{itemize}
        \item The definition is independent of the choice of the good resolution. In particular, these invariants are intrinsic to $X$.
        \item If $X$ has only terminal quotient singularities, then the definition is compatible with the definition using smooth orbifolds in \cite{Self_LG_Model} and \cite{Self_LG_divisorial}.
        \item If $X$ has a $\mathbb{Q}$-smoothing, then the definition is compatible with the definition using $\mathbb{Q}$-smoothings in \cite{Self_LG_Model} and \cite{Self_LG_divisorial}.
    \end{itemize}
    In particular, we generalize the definition of the quantum period to general terminal Fano threefolds.

\subsection{Definition}\label{subsec:stationary-definition}

    We first specify the resolutions and the stationary data used in the construction.

\begin{definition}[Good resolutions]\label{def:stationary-good-resolution}
    Let $X$ be a projective terminal Fano threefold, and put
    $$
    X_{sm}=X\setminus\mathrm{Sing}(X).
    $$
    A \emph{good resolution} of $X$ is a projective birational morphism $p:Y\rightarrow X$ such that $Y$ is smooth, $p$ is an isomorphism over $X_{sm}$, and
    $$
    E=\left(p^{-1} \mathrm{Sing}(X) \right)_{\mathrm{red}} =\sum\limits_{i=1}^{r}E_{i}
    $$
    is a simple-normal-crossings divisor.
\end{definition}

\begin{definition}[Stationary descendants on a good resolution]\label{def:stationary-resolution-invariant}
    Notation as in \cref{def:stationary-good-resolution}. We use connected stable maps and rational cohomology. Choose compatible canonical divisors and write
    $$
    K_{Y}=p^{*}K_{X}+\sum\limits_{i=1}^{r}a_{i}E_{i}
    $$
    where $E_{i}$ are exceptional divisors over $X$ and $a_{i}>0$. Fix integers $g\geq 0$ and $n\geq 1$, homogeneous classes $\gamma_{j}\in H^{2c_{j}}(X,\mathbb{Q})$, $c_{j}\in\mathbb{Z}_{>0}$, $1 \leq j\leq n$, and nonnegative integers $k_{j}$. Let
    $$
    \beta\in H_{2}(X,\mathbb{Z}), \qquad d=(-K_{X})\cdot\beta.
    $$
    We consider data satisfying the dimension condition
    $$
    \sum\limits_{j=1}^{n}(k_{j}+c_{j})=d+n.
    $$
    In particular, $d$ is a nonnegative integer. We use the usual stable range: when $\beta=0$, we require $2g-2+n>0$. 
    
    Define
    $$
    \begin{aligned}
    \left\langle
    \prod_{j=1}^{n}\tau_{k_{j}}(\gamma_{j})
    \right\rangle^{X,p}_{g,n,\beta}
    :={}&
    \sum\limits_{\substack{
    \alpha\in H_{2}(Y,\mathbb{Z})\\
    p_{*}\alpha=\beta\\
    E_{i}\cdot\alpha=0\ \text{for every }i}}
    \left\langle
    \prod_{j=1}^{n}\tau_{k_{j}}(p^{*}\gamma_{j})
    \mid\emptyset
    \right\rangle^{(Y,E),\log}_{g,n,\alpha}.
    \end{aligned}
    $$
    The relative invariant on the right has $n$ absolute markings, all with contact order zero, and no relative markings. Its descendants are pulled back from the underlying unexpanded stable-map space.
    
    Classes with no representing stable maps contribute zero. Finiteness of the nonzero summands and independence of $p$ are proved in \cref{thm:stationary-resolution-independence}.
    
    If the dimension condition fails, the invariant is defined to be zero. The invariant on the right-hand side is that of the full logarithmic relative moduli space, including its boundary strata.
\end{definition}

\subsection{Two log smooth degenerations}\label{subsec:stationary-degenerations}

    The following constructions provide the two degenerations used in the blow-up comparison.

\begin{construction}[The first log smooth degeneration]\label{construction: deformation family 1}
    Let $Z$ be a smooth projective threefold, $B=\sum\limits_{i=1}^{r}B_{i}$ a reduced simple normal crossings divisor, and $A$ a nef $\mathbb{Q}$-Cartier divisor on $Z$. Put $T=\mathbb{A}^{1}$, with coordinate $t$. The total space and the base carry the divisorial logarithmic structures of the central fibre and the origin. We construct a projective log smooth degeneration
    $$
    \pi:\mathcal{W}\rightarrow T
    $$
    and a relatively nef $\mathbb{Q}$-Cartier divisor $\mathcal{L}$ with the following properties.
    
    \begin{enumerate}
    \item The total space $\mathcal{W}$ is a smooth scheme, and $\mathcal{W}_{0}$ is a reduced simple normal crossings divisor.
    
    \item For $t\neq0$, the fibre $\mathcal{W}_t$ is isomorphic to $Z$.
    
    \item The central fibre contains a distinguished component $Z_{0}\simeq Z$ whose reduced vertical boundary is $B$. Thus, under this identification,
    $$
    \left( Z_{0}\cap\overline{\mathcal{W}_{0}\setminus Z_{0}} \right)_{\mathrm{red}} = B.
    $$
    
    \item There is a projective birational morphism
    $$
    \rho:\mathcal{W}\rightarrow Z\times T
    $$
    over $T$, which is an isomorphism away from $\mathrm{Supp}(B)\times\{0\}$. Writing $h=\mathrm{pr}_{Z}\circ\rho$, we have
    $$
    h|_{Z_{0}}=\mathrm{id}_{Z}.
    $$
    Every other central component maps into $\mathrm{Supp}(B)$ under $h$.
    
    \item
    The divisor $\mathcal{L}=h^{*}A$ satisfies
    $$
    \mathcal{L}|_{\mathcal{W}_t}=A \quad(t\neq0), \qquad \mathcal{L}|_{Z_{0}}=A.
    $$
    \end{enumerate}

    Indeed, start with
    $$
    \mathcal{W}^{(0)}=Z\times T, \qquad Z^{(0)}=Z\times\{0\}.
    $$
    For $i=1,\ldots,r$, let $Z^{(i-1)}$ be the strict transform of $Z^{(0)}$ and let $H_{i}^{(i-1)}$ be the strict transform of $B_{i}\times T$ on $\mathcal{W}^{(i-1)}$. Define
    $$
    \mathcal{W}^{(i)} = \mathrm{Bl}_{Z^{(i-1)}\cap H_{i}^{(i-1)}} \mathcal{W}^{(i-1)}.
    $$
    
    At every step the center is a smooth divisor on the distinguished component. Hence the strict transformation of $Z$ on the blow-up remains isomorphic to $Z$, while $B_{i}$ is added to its vertical boundary. The central fibre and the remaining horizontal divisors continue to have simple normal crossings. The horizontal divisors already separated from the distinguished component remain disjoint from it.
    
    Set $\mathcal{W}=\mathcal{W}^{(r)}$, and let $\rho$ be the composite blow-up morphism. All centers lie over $\mathrm{Supp}(B)\times\{0\}$, so the general fibres are unchanged and all auxiliary central components map into $\mathrm{Supp}(B)$. The reduced vertical boundary on $Z_{0}$ is exactly $B$.
    
    Finally, the pullback $h^{*}A$ is relatively nef and has the stated restrictions.
\end{construction}

\begin{construction}[The second log smooth degeneration]\label{construction: deformation family 2}
    Notation as in \cref{construction: deformation family 1}. Let
    $$
    f:\widetilde{Z}=\mathrm{Bl}_C Z\rightarrow Z
    $$
    be the blow-up of a smooth closed center of codimension at least two having normal crossings with $B$. Let $P=\mathrm{Exc}(f)_{\mathrm{red}}$, and put
    $$
    \widetilde{B}=P+\sum\limits_{i=1}^{r}\widetilde{B}_{i},
    $$
    where $\widetilde{B}_{i}$ is the strict transform of $B_{i}$.

    We construct a projective log smooth degeneration
    $$
    \pi':\mathcal{W}'\rightarrow T
    $$
    and a relatively nef $\mathbb{Q}$-Cartier divisor $\mathcal{L}'$ with the following properties.

    \begin{enumerate}
    \item The total space $\mathcal{W}'$ is a smooth scheme, and $\mathcal{W}'_{0}$ is a reduced simple normal crossings divisor.

    \item For $t\neq0$, the fibre $\mathcal{W}'_t$ is isomorphic to $Z$.

    \item The central fibre contains a distinguished component $Z_{0}^{+}\simeq \widetilde{Z}$ whose reduced vertical boundary is $\widetilde{B}$. Thus, under this identification,
    $$
    \left( Z_{0}^{+}\cap \overline{\mathcal{W}'_{0}\setminus Z_{0}^{+}} \right)_{\mathrm{red}} = \widetilde{B}.
    $$

    \item There is a projective birational morphism
    $$
    \rho':\mathcal{W}'\rightarrow Z\times T
    $$
    over $T$, which is an isomorphism away from $(\mathrm{Supp}(B)\cup C)\times\{0\}$. Writing $h'=\mathrm{pr}_{Z}\circ\rho'$, we have
    $$
    h'|_{Z_{0}^{+}}=f.
    $$
    Every other central component maps into $\mathrm{Supp}(B)\cup C$ under $h'$.

    \item The divisor $\mathcal{L}'=(h')^{*}A$ satisfies
    $$
    \mathcal{L}'|_{\mathcal{W}'_t}=A \quad(t\neq0), \qquad \mathcal{L}'|_{Z_{0}^{+}}=f^{*}A.
    $$
    \end{enumerate}

    Indeed, begin with the deformation to the normal cone
    $$
    \mathcal{V} = \mathrm{Bl}_{C\times\{0\}}(Z\times T) \rightarrow T.
    $$
    The strict transform of $Z\times\{0\}$ is $\widetilde{Z}$. If $Q$ is the exceptional divisor of this blow-up, then
    $$
    \mathcal{V}_{0}=\widetilde{Z}\cup_P Q.
    $$
    The morphism to $Z$ restricts to $f$ on $\widetilde{Z}$, and $Q$ maps into $C$.
    
    To verify semistability, choose local parameters $x_{1},x_{2},x_3$ with
    $$
    \mathcal{I}_C=(x_{1},\ldots,x_{c}), \qquad c=\mathrm{codim}_{Z} C.
    $$
    The blow-up ideal on $Z\times T$ is $(t,x_{1},\ldots,x_{c})$. In the $x_{a}$-chart, for $1\leq a\leq c$, write
    $$
    t=x_{a} u, \qquad x_{j}=x_{a} y_{j} \quad (1\leq j\leq c,\ j\neq a).
    $$
    The central fibre is $x_{a} u=0$, with distinguished component $u=0$ and exceptional component $x_{a}=0$. In the $t$-chart, write $x_{j}=t y_{j}$. The central fibre is the smooth divisor $t=0$.
    
    Since $C$ has normal crossings with $B$, these coordinates also show that the central fibre and the strict transforms of $B_{i}\times T$ have simple normal crossings.
    
    The restriction of the strict transform of $B_{i}\times T$ to $\widetilde{Z}$ is $\widetilde{B}_{i}$. Successively blow up the current total space along the intersection of its distinguished component with the strict transform of $B_{i}\times T$, as in \cref{construction: deformation family 1}. The distinguished component remains isomorphic to $\widetilde{Z}$, and its reduced vertical boundary becomes exactly
    $$
    P+\sum\limits_{i} \widetilde{B}_{i}=\widetilde{B}.
    $$
    
    Let $\mathcal{W}'$ be the resulting total space and $\rho'$ its morphism to $Z\times T$. All further centers lie over $\mathrm{Supp}(B)\times\{0\}$. This proves the assertions about the fibres and the images of the auxiliary components. The claims concerning $\mathcal{L}'=(h')^{*}A$ follow from nefness of $A$ and the restrictions of $h'$.
    
    If $C\subseteq\mathrm{Supp}(B)$, then
    $$
    \widetilde{B}=(f^{-1}B)_{\mathrm{red}}.
    $$
\end{construction}

\subsection{Independence of the good resolution}\label{subsec:stationary-resolution-independence}

    In this subsection, we show that the above definition is independent of the choice of the good resolution. We first explain the role of the positive-degree insertion condition.

\begin{lemma}[Support of the insertions]\label{lem:stationary-support}
    Let $X$ be a projective terminal Fano threefold, and put $X_{sm}=X\setminus\mathrm{Sing}(X)$. Then for every integer $c>0$, the natural map
    $$
    H_{c}^{2c}(X_{sm},\mathbb{Q}) \rightarrow H^{2c}(X,\mathbb{Q})
    $$
    is an isomorphism.
    
    Moreover, let $p:Y\rightarrow X$ be a good resolution, and put
    $$
    E=(p^{-1}\mathrm{Sing}(X))_{\mathrm{red}}.
    $$
    Then for every $\gamma\in H^{2c}(X,\mathbb{Q})$, the class $p^{*}\gamma$ is represented by the corresponding compactly supported class on $Y\setminus E\simeq X_{sm}$. Its restriction to $E$ is zero. More generally, if a morphism $a:A\rightarrow X$ factors through $\mathrm{Sing}(X)$, then $a^{*}\gamma=0$.
\end{lemma}

\begin{proof}
    Since $X$ is proper, the localization sequence gives
    $$
    H^{2c-1}(\mathrm{Sing}(X),\mathbb{Q}) \rightarrow H_{c}^{2c}(X_{sm},\mathbb{Q}) \rightarrow H^{2c}(X,\mathbb{Q}) \rightarrow H^{2c}(\mathrm{Sing}(X),\mathbb{Q}).
    $$
    Both outer groups vanish because $\mathrm{Sing}(X)$ is finite and $c>0$. The assertion about pullback follows from the naturality of the localization sequence and the identification $Y\setminus E\simeq X_{sm}$. The remaining assertions follow because a positive-degree class restricts to zero on the finite singular locus.
\end{proof}

    We next extend the marked-factor vanishing argument to several markings and arbitrary genus.

\begin{lemma}[Stationary degeneration comparison]\label{lem:stationary-degeneration}
    Let $X$ be a projective terminal Fano threefold, and let $\mathcal{W}\rightarrow T$ be a projective semistable degeneration over a smooth pointed curve $0\in T$.

    Let $Z$ be its distinguished component, and write its reduced vertical boundary as
    $$
    B= \left( Z\cap\overline{\mathcal{W}_{0}\setminus Z} \right)_{\mathrm{red}} = \sum\limits_{i} B_{i}.
    $$
    Let $p:Z\rightarrow X$ be projective and birational, isomorphic over $X_{sm}=X\setminus\mathrm{Sing}(X)$. Assume that the $B_{i}$ are exactly the exceptional prime divisors of $p$, so that
    $$
    K_{Z}=p^{*}K_X+\sum\limits_{i} a_{i}B_{i}, \qquad a_{i}>0.
    $$

    Let $\mathcal{L}$ be a relatively nef $\mathbb{Q}$-Cartier divisor on $\mathcal{W}$ with
    $$
    \mathcal{L}|_{Z}=p^{*}(-K_X).
    $$
    Fix $g\geq 0$, $n\geq1$, homogeneous classes $\gamma_{j}\in H^{2c_{j}}(X,\mathbb{Q})$ with integers $c_{j}>0$, and integers $k_{j}\geq 0$. Assume that the insertion classes extend over the degeneration, restrict to $p^{*}\gamma_{j}$ on $Z$, and restrict to zero on every other central component. The descendants use the coarse cotangent-line classes pulled back from the underlying unexpanded stable-map space.
    
    For connected stable maps of genus $g$ and total $\mathcal{L}$-degree $d$ satisfying
    $$
    \sum\limits_{j=1}^{n}(k_{j}+c_{j})=d+n,
    $$
    all contributions with relative marked points on the distinguished factors vanish. The remaining contribution is the sum of logarithmic invariants of $(Z,B)$ with no relative marked points, with coefficient one.
\end{lemma}

\begin{proof}
    A contribution placing an original marking on an auxiliary component vanishes by the restriction of its insertion. Thus all original insertions occur on factors assigned to $Z$.
    
    Let $\mathcal{V}_{Z}$ be the collection of these connected relative factors, including the compatible refinements used for cutting and gluing. For $v\in\mathcal{V}_{Z}$, write $\alpha_{v}$ for its class, $n_{v}$ for its number of original markings, and $l_{v}$ for its number of relative markings.
    
    For the factor indexed by $v$, let $(Z_v,B_v)$ be its refined distinguished target pair, let $p_v:Z_v\rightarrow X$ be the induced morphism, and write $B_v=\sum\limits_{i} B_{v,i}$. Let $\mathcal{M}_v$ be its relative moduli stack. Write
    $$
    K_{Z_{v}}=p_{v}^{*}K_{X}+\sum\limits_{i} a_{v,i}B_{v,i}.
    $$
    Every boundary component is exceptional over $X$, so $a_{v,i}>0$. If the contact orders are $\mu_{h,i}$, then \cref{prop:relative-virtual-dimension} gives
    $$
    \mathrm{vdim}\mathcal{M}_{v} = p_{v}^{*}(-K_{X})\cdot\alpha_{v}+n_{v} + \sum\limits_{h=1}^{l_{v}} \left( 1-\sum\limits_{i}(a_{v,i}+1)\mu_{h,i} \right).
    $$
    There is no genus term in dimension three. Every term in parentheses is strictly negative.
    
    Let $\delta\geq 0$ be the complementary $\mathcal{L}$-degree. Then
    $$
    \sum\limits_{v\in\mathcal{V}_{Z}} p_{v}^{*}(-K_{X})\cdot\alpha_{v}=d-\delta, \qquad \sum\limits_{v\in\mathcal{V}_{Z}}n_{v}=n.
    $$
    If any relative marking occurs, then
    $$
    \sum\limits_{v\in\mathcal{V}_{Z}} \mathrm{vdim}\mathcal{M}_{v}<d+n.
    $$
    The total insertion degree is $2(d+n)$. Hence at least one distinguished factor receives insertions whose degree exceeds twice its virtual complex dimension.

    By \cref{lem:stationary-support}, the original insertions can be represented away from the exceptional boundary. Their descendant classes therefore pull back from the indicated factors. \cref{lem:scheme-factor-vanishing} eliminates the glued contribution.

    When no relative marking occurs, \cref{lem:ordinary-edge-free} identifies the remaining contribution.
\end{proof}

\begin{lemma}[Comparison under a blow-up]\label{lem:stationary-blowup}
    Notation as in \cref{def:stationary-resolution-invariant}. Let
    $$
    f:\widetilde{Y}=\mathrm{Bl}_{C}Y\rightarrow Y
    $$
    be the blow-up of a smooth connected center $C\subset E$ of codimension at least two, having normal crossings with $E$. Write $\widetilde{E}=(f^{-1}E)_{\mathrm{red}} =\sum\limits_{j}\widetilde{E}_{j}$, including the exceptional divisor of $f$ among these components. Then, for any $\alpha\in H_{2}(Y,\mathbb{Z})$ satisfying $p_{*}\alpha=\beta$ and $E_{i}\cdot\alpha=0$ for every $i$, we have
    $$
    \left\langle
    \prod_{j=1}^{n}\tau_{k_{j}}(p^{*}\gamma_{j})
    \mid\emptyset
    \right\rangle^{(Y,E),\log}_{g,n,\alpha} =
    \sum\limits_{\substack{f_{*}\widetilde{\alpha}=\alpha\\
    \widetilde{E}_{j}\cdot\widetilde{\alpha}=0\ \text{for every }j}}
    \left\langle\prod_{j=1}^{n}\tau_{k_{j}}((p \circ f)^{*}\gamma_{j})\mid\emptyset
    \right\rangle^{(\widetilde{Y},\widetilde{E}),\log}_{g,n,\widetilde{\alpha}},
    $$
    where $\widetilde{\alpha}\in H_{2}(\widetilde{Y},\mathbb{Z})$ and $f_{*}\widetilde{\alpha}=\alpha$ is imposed in integral homology.
\end{lemma}

\begin{proof}
    For a class in the displayed sum, the boundary intersection conditions and the relative canonical formula give
    $$
    (-K_Y-E)\cdot\alpha = (-K_Y)\cdot\alpha = (-K_X)\cdot\beta = d.
    $$
    By \cref{prop:relative-virtual-dimension}, the virtual complex dimension is $d+n$, as required by the insertion-degree condition. Apply \cref{construction: deformation family 1,construction: deformation family 2} to $(Y,E)$ and $f:\widetilde{Y}\rightarrow Y$, with $A=p^{*}(-K_X)$. The first degeneration has general fibre $Y$ and distinguished pair $(Y,E)$. The second starts with deformation to the normal cone of $C$ and has general fibre $Y$ and distinguished pair $(\widetilde{Y},\widetilde{E})$. Both families have a morphism to $Y$, and hence to $X$. Every auxiliary component maps into $\mathrm{Sing}(X)$, since all blow-up centers lie over $\mathrm{Sing}(X)$. Equip the families with the pullback of $p^{*}(-K_{X})$, and extend the insertions by pullback from $X$.
    
    \cref{lem:stationary-degeneration} applied to the first degeneration gives
    $$
    \left\langle \prod_{j=1}^{n}\tau_{k_{j}}(p^{*}\gamma_{j}) \right\rangle^{Y}_{g,n,\alpha} =
    \left\langle \prod_{j=1}^{n}\tau_{k_{j}}(p^{*}\gamma_{j}) \mid\emptyset \right\rangle^{(Y,E),\log}_{g,n,\alpha}.
    $$
    Applied to the second degeneration, it gives
    $$
    \left\langle \prod_{j=1}^{n}\tau_{k_{j}}(p^{*}\gamma_{j}) \right\rangle^{Y}_{g,n,\alpha}= 
    \sum\limits_{\substack{f_{*}\widetilde{\alpha}=\alpha\\
    \widetilde{E}_{j}\cdot\widetilde{\alpha}=0\ \text{for every }j}}
    \left\langle\prod_{j=1}^{n}\tau_{k_{j}}((p \circ f)^{*}\gamma_{j})\mid\emptyset
    \right\rangle^{(\widetilde{Y},\widetilde{E}),\log}_{g,n,\widetilde{\alpha}}.
    $$
    Combining these identities proves the assertion.

    Both total families map to $Y$, so we retain $\alpha\in H_{2}(Y,\mathbb{Z})$ itself. The integral pushforward of a stable map to the fixed target is locally constant in families. On each bounded mapping problem, fixing this class selects open and closed substacks, and the virtual identities restrict to them. Thus the surviving classes satisfy $f_{*}\widetilde{\alpha}=\alpha$ in integral homology (cf. \cite[Theorem~3.3 and Remark~3.4]{LiuYau04}).

    Every $f^{*}E_{i}$ is supported on $\widetilde{E}$. Hence $\widetilde{E}_{j}\cdot\widetilde{\alpha}=0$ for every $j$ implies $E_{i}\cdot f_{*}\widetilde{\alpha} =f^{*}E_{i}\cdot\widetilde{\alpha}=0$ for every $i$.

    The comparisons are finite for fixed $\alpha$. Indeed, let $P$ be the exceptional divisor of $f$. An ample divisor on $\widetilde{Y}$ differs from the pullback of a divisor on the smooth variety $Y$ by a multiple of $P$. Its degree is therefore fixed on classes satisfying
    $$
    f_{*} \widetilde{\alpha}=\alpha, \qquad P\cdot \widetilde{\alpha}=0.
    $$
    Boundedness of stable maps gives finiteness.

    If $l$ is a line in a fibre of $P\rightarrow C$, then
    $$
    \ker\!\left( f_{*}:H_{2}(\widetilde{Y},\mathbb{Z}) \rightarrow H_{2}(Y,\mathbb{Z}) \right) = \mathbb{Z}[l], \qquad P\cdot l=-1.
    $$
    Consequently, the integral pushforward and zero degree against $P$ determine at most one integral lift.
\end{proof}

\begin{theorem}[Resolution independence]\label{thm:stationary-resolution-independence}
    Notation as in \cref{def:stationary-resolution-invariant}. Then the sum defining
    $$
    \left\langle \prod_{j=1}^{n}\tau_{k_j}(\gamma_j) \right\rangle^{X,p}_{g,n,\beta}
    $$
    has finitely many nonzero summands, and its value is independent of the good resolution $p$.
    
    Moreover, for fixed $X$, $g$, $n$, $d$, and $p$, only finitely many integral classes on $Y$ can contribute as $\beta$, the stationary insertions, and the descendant exponents vary subject to these fixed values and the dimension condition.
\end{theorem}

\begin{proof}
    We first establish a degree bound on one appropriate good resolution. Put $X_{sm}=X\setminus\mathrm{Sing}(X)$. By \cite[Theorem 1, Remark 5]{KW21}, there is a projective good resolution $p_{0}:Y_{0}\rightarrow X$ whose exceptional locus supports a $p_{0}$-ample divisor
    $$
    F=\sum\limits_{i} f_{i}E_{0,i}.
    $$
    Write $E_{0}=\sum\limits_{i}E_{0,i}$ for its reduced exceptional divisor. The divisor $F$ is not required to be effective. Choose a sufficiently large divisible integer $m$ such that
    $$
    A_{Y_{0}}=m\,p_{0}^{*}(-K_{X})+F
    $$
    is ample and Cartier. If $E_{0,i}\cdot\alpha_{0}=0$ for all $i$, and $(-K_{X})\cdot(p_{0})_{*}\alpha_{0}=d$, then
    $$
    A_{Y_{0}}\cdot\alpha_{0}=md.
    $$
    Stable maps of the fixed genus and number of markings and bounded $A_{Y_{0}}$-degree form a finite-type stack. Their integral homology class is locally constant, so only finitely many integral classes $\alpha_{0}$ occur. Hence the required sums are finite on $Y_{0}$.

    Now let $p:Y\rightarrow X$ be any good resolution. We apply the weak factorization to the birational map $Y\dashrightarrow Y_{0}$ over the common open set $X_{sm}$ (cf. \cite{AKMW02}). The intermediate varieties are smooth and projective, the centers are disjoint from $X_{sm}$ and have normal crossings with the respective boundaries, and every intermediate variety has a projective morphism to $Y$ or $Y_{0}$. In particular, every intermediate variety has a projective morphism to $X$. Hence every intermediate pair is again a good resolution of $X$. Without loss of generality, we can assume that every smooth center is connected. We can then apply \cref{lem:stationary-blowup} to every nontrivial step of the factorization. Each elementary comparison retains the integral pushforward to $X$. The uniqueness of the integral lift with zero new exceptional degree, proved in \cref{lem:stationary-blowup}, transports finite support through the factorization. Summing the coefficient comparisons over all classes with pushforward $\beta$ gives
    $$
    \left\langle \prod_{j=1}^{n}\tau_{k_{j}}(\gamma_{j}) \right\rangle^{X,p}_{g,n,\beta} =
    \left\langle \prod_{j=1}^{n}\tau_{k_{j}}(\gamma_{j}) \right\rangle^{X,p_{0}}_{g,n,\beta}.
    $$
    This proves finiteness and independence.
    
    For the last assertion, choose homogeneous bases of the finite-dimensional positive even cohomology of $X$. For fixed $n$ and $d$, there are only finitely many choices of their degrees and of the nonnegative descendant exponents satisfying the dimension condition. Taking the finite union of the supports for these basis insertions proves the assertion.
\end{proof}

\begin{convention}\label{conv:resolution-invariants}
    In the rest of the paper, we write
    $$
    \left\langle \prod_{j=1}^{n}\tau_{k_{j}}(\gamma_{j}) \right\rangle^{X,\mathrm{res}}_{g,n,\beta}
    $$
    for the invariant in \cref{def:stationary-resolution-invariant}, omitting the resolution.
\end{convention}

\subsection{Compatibility with orbifold Gromov-Witten invariants}\label{subsec:quotient-resolution-comparison}

    Let $X$ be a projective terminal quotient Fano threefold, and let $\mathfrak{X}$ be its canonical covering stack. There are two different ways to define its Gromov-Witten invariants: using \cref{definition: Gromov-Witten for orbifold} and \cref{def:stationary-resolution-invariant}, respectively. The main goal of this subsection is to show that these two definitions coincide.

\begin{proposition}[Comparison with orbifold invariants]\label{prop:quotient-resolution-comparison}
    Let $X$ and the stationary data $g,n,\beta,d,\gamma_{j},c_{j},k_{j}$ be as in \cref{thm:stationary-resolution-independence}. Assume in addition that $X$ has terminal quotient singularities, and let $\mathfrak{X}$ be its canonical covering stack. Regard the $\gamma_{j}$ as untwisted insertions on $\mathfrak{X}$. For the fixed integral class $\beta$, we have
    $$
    \left\langle \prod_{j=1}^{n}\tau_{k_{j}}(\gamma_{j}) \right\rangle^{X,\mathrm{res}}_{g,n,\beta}
    =\left\langle \prod_{j=1}^{n}\tau_{k_{j}}(\gamma_{j}) \right\rangle^{\mathfrak{X}}_{g,n,\beta}.
    $$
\end{proposition}

\begin{proof}
    We first consider a single Kawamata blow-up. Let $V$ be a projective terminal quotient threefold, let $P\in V$ be a terminal quotient singularity of type $\frac{1}{r}(a,r-a,1)$ for some positive integers $a<r$ where $\gcd(a,r)=1$, and let
    $$
    f:\widetilde{V}\rightarrow V
    $$
    be its Kawamata blow-up, with exceptional divisor $E$. Then
    $$
    K_{\widetilde{V}}=f^{*}K_{V}+\frac{1}{r}E.
    $$
    We apply \cite[Construction~3.6]{Self_LG_divisorial} to this divisorial contraction. More explicitly, we take the weighted blow-up of the trivial family at $(P,0)$ with weight $\frac{1}{r}(a,r-a,1,r)$. Its general fibre is the canonical covering stack $\mathfrak{V}$. Its distinguished pair is $(\widetilde{\mathfrak{V}},\mathfrak{E})$, where $\widetilde{\mathfrak{V}}$ is the canonical covering stack of $\widetilde{V}$ and $\mathfrak{E}$ is a smooth divisor. The auxiliary component maps to $P$.

    Now we apply \cref{theorem: orbifold degeneration formula} to this family. Take $n\geq1$, integers $k_{i}\geq 0$, and homogeneous untwisted insertions
    $$
    \gamma_{i}\in H^{2c_{i}}(V,\mathbb{Q}), \qquad c_{i}\in\mathbb{Z}_{>0},
    $$
    represented away from $P$. Let $\beta\in H_{2}(V,\mathbb{Z})$. Assume
    $$
    \sum\limits_{i=1}^{n}(k_{i}+c_{i})=d+n, \qquad d=(-K_{V})\cdot\beta.
    $$
    Fix a splitting
    $$
    \eta=(\Xi_{1},\Xi_{2})\in\Omega.
    $$
    The sets $N_{1},N_{2}$ index the original markings, $M$ labels the relative markings, and $d_{j}>0$ is the intersection multiplicity at $j\in M$. A term with $N_{2}\neq\emptyset$ vanishes because the auxiliary component maps to $P$. Hence we may assume
    $$
    N_{1}=\{1,\ldots,n\}, \qquad N_{2}=\emptyset.
    $$
    Let $\beta_{1}$ be the total coarse curve class specified by $\Xi_{1}$, and let $\zeta_{j}$ be the monodromy at $j\in M$. Then we have
    $$
    f^{*}(-K_{V})\cdot\beta_{1}=d, \qquad E\cdot\beta_{1}=\sum\limits_{j\in M}d_{j}.
    $$

    The normal monodromy has fractional weight $\{d_{j}\}$. Hence
    $$
    \mathrm{Age} \left( \zeta_{j}\mid T_{\widetilde{\mathfrak{V}}}(-\log\mathfrak{E}) \right) = \mathrm{Age}_{\widetilde{\mathfrak{V}}}(\zeta_{j}) -\{d_{j}\}.
    $$
    Fix a twisting choice $\mathfrak{r}$. Then by the orbifold Riemann-Roch formula, we have
    $$
    \mathrm{vdim} \mathcal{K}_{\Xi_{1}}^{\mathfrak{r}} (\widetilde{\mathfrak{V}},\mathfrak{E}) =
    d+n+\sum\limits_{j\in M}
    \left( 1-\mathrm{Age}_{\widetilde{\mathfrak{V}}}(\zeta_{j}) -\frac{d_{j}}{r}-\lfloor d_{j}\rfloor \right).
    $$
    For every $j \in M$, we have
    $$
    \mathrm{Age}_{\widetilde{\mathfrak{V}}}(\zeta_{j}) +\frac{d_{j}}{r}+\lfloor d_{j}\rfloor>1.
    $$
    Indeed, for nontrivial $\zeta_{j}$, this follows from \cref{thm:reid-tai}. For trivial $\zeta_{j}$, the intersection multiplicity $d_{j}$ is a positive integer and the expression is $(1+1/r)d_{j}>1$.
    
    If $M \neq \emptyset$, the distinguished-side virtual dimension is therefore less than $d+n$. The original descendants already have complex codimension $d+n$, while the relative insertions have nonnegative ordinary cohomological degree. Thus the distinguished-side relative invariant, and hence the corresponding degeneration term, is zero. Hence only $M=\emptyset$ remains. Since the original source is connected and carries an original marking, the surviving source is connected and lies on the distinguished side. The gluing coefficient is one.

    The coarse total space of the weighted degeneration maps to $V$. Retain the exact integral pushforward $\beta\in H_{2}(V,\mathbb{Z})$ in deformation invariance and in the degeneration formula, by the same open-and-closed restriction used in \cref{lem:stationary-blowup}. Thus the surviving classes satisfy $f_{*}\alpha=\beta$ in integral homology, and we obtain
    $$
    \left\langle \prod_{j}\tau_{k_{j}}(\gamma_{j}) \right\rangle^{\mathfrak{V}}_{g,n,\beta} =
    \sum\limits_{\substack{ \alpha\in H_{2}(\widetilde{V},\mathbb{Z})\\ f_{*}\alpha=\beta\\ E\cdot\alpha=0}}
    \left\langle \prod_{j}\tau_{k_{j}}(f^{*}\gamma_{j}) \mid\emptyset \right\rangle^{(\widetilde{\mathfrak{V}},\mathfrak{E})}_{g,n,\alpha}.
    $$

    Next apply deformation to the normal cone of the smooth divisor $\mathfrak{E}\subset\widetilde{\mathfrak{V}}$. Keeping $\alpha\in H_{2}(\widetilde{V},\mathbb{Z})$ fixed through the morphism to $\widetilde{V}$, for $E\cdot\alpha=0$, the same calculation, using $f^{*}(-K_{V})$ for degree conservation, eliminates every nonempty relative datum. Hence
    $$
    \left\langle\prod_{j}\tau_{k_{j}}(f^{*}\gamma_{j})\mid\emptyset \right\rangle^{(\widetilde{\mathfrak{V}},\mathfrak{E})}_{g,n,\alpha} =
    \left\langle \prod_{j}\tau_{k_{j}}(f^{*}\gamma_{j}) \right\rangle^{\widetilde{\mathfrak{V}}}_{g,n,\alpha}.
    $$

    We take the good resolution $p: Y\rightarrow X$ of $X$ to be the successive Kawamata blow-ups on the remaining singular points. At each step, the new exceptional divisor has degree zero. The zero-degree conditions for previous exceptional divisors are preserved by their total-transform formulas. The insertions are pulled back from $X$ and remain supported away from all extraction centers. Write
    $$
    E=(p^{-1}\mathrm{Sing}(X))_{\mathrm{red}} =\sum\limits_{i}E_{i}.
    $$
    The successive comparisons give
    $$
    \left\langle \prod_{j}\tau_{k_{j}}(\gamma_{j}) \right\rangle^{\mathfrak{X}}_{g,n,\beta} =
    \sum\limits_{\substack{ \alpha\in H_{2}(Y,\mathbb{Z})\\ p_{*}\alpha=\beta\\ E_{i}\cdot\alpha=0\ \forall i}} 
    \left\langle \prod_{j}\tau_{k_{j}}(p^{*}\gamma_{j}) \right\rangle^{Y}_{g,n,\alpha}.
    $$
    Apply \cref{construction: deformation family 1} to $(Y,E)$ with $A=p^{*}(-K_X)$. Together with \cref{lem:stationary-degeneration}, it identifies this sum with the logarithmic sum defining the resolution invariant.
    
    For a single Kawamata blow-up, an ample divisor on $\widetilde{V}$ differs from the pullback of a $\mathbb{Q}$-Cartier divisor on $V$ by a multiple of $E$. The conditions $f_{*}\alpha=\beta$ and $E\cdot\alpha=0$ therefore fix its degree. Boundedness and local constancy of the integral class give finitely many integral classes in each comparison. The final sum is finite by \cref{thm:stationary-resolution-independence}.
\end{proof}

\subsection{Compatibility with a \texorpdfstring{$\mathbb{Q}$}{Q}-smoothing}\label{subsec:stationary-smoothing}

    We recall the definition of $\mathbb{Q}$-smoothing.

\begin{definition}[$\mathbb{Q}$-smoothings] \label{def:Q-smoothing}
    Let $X$ be a normal $\mathbb{Q}$-Gorenstein variety. A \emph{$\mathbb{Q}$-smoothing} of $X$ is a flat $\mathbb{Q}$-Gorenstein family
    $$
    \pi:\mathcal{X}\rightarrow C,\qquad \mathcal{X}_{0}\simeq X,
    $$
    over a smooth pointed curve $0\in C$, such that $\mathcal{X}_{s}$ has only quotient singularities for $s\neq 0$, after shrinking $C$ around $0$.
    
    We also use the analytic version of this definition, with $C$ replaced by a sufficiently small disc. When $X$ is proper, the family is required to be proper. We call $X$ \emph{$\mathbb{Q}$-smoothable} if it admits a $\mathbb{Q}$-smoothing.
    
    In the terminal threefold setting, we shrink the base so that the nearby fibres are terminal. Their quotient singularities are therefore isolated, as required in \cite[Definition~1.2]{Sano15}.
\end{definition}
    
   We first extend the line bundles needed for the divisor weights.
    
\begin{lemma}[Extension of divisor classes]\label{lem:extend-divisor-classes}
    Let $X$ be a projective terminal Fano threefold, and let
    $$
    \pi:\mathcal{X}\rightarrow T \ni 0, \qquad \mathcal{X}_{0}\simeq X,
    $$
    be a projective $\mathbb{Q}$-smoothing over a smooth pointed curve $0\in T$. Then any prescribed finite collection of line bundles on $X$ extends to $\mathcal{X}$ after shrinking $T$ and making an \'etale base change.
\end{lemma}

\begin{proof}
    Since $X$ is a terminal Fano threefold, we have
    $$
    H^{1}(X,\mathcal{O}_{X}) = H^{2}(X,\mathcal{O}_{X}) = 0.
    $$
    Hence the relative Picard functor extends each chosen class \'etale locally on $T$. After a further \'etale localization, the classes are represented by line bundles. A common base change works for the finite collection.
\end{proof}

    The next lemma provides a semistable model for the comparison of the central and general invariants.

\begin{lemma}[A semistable reduction]\label{lem:scheme-resolution-family}
    Notation as in \cref{lem:extend-divisor-classes}. Then after shrinking $T$ and taking a finite base change, there is a projective morphism
    $$
    p:\mathcal{Y}\rightarrow\mathcal{X}
    $$
    over the new base. Here and below in this lemma, $T$ and $\mathcal{X}$ denote the new base and the base-changed family, and $\mathcal{Y}_{s}$ is the fibre of $\pi\circ p$ over $s$. The morphism has the following properties:
    \begin{enumerate}
    \item The family $\mathcal{Y}/T$ is semistable.
    \item For $s\neq 0$, the induced morphism $p_{s}:\mathcal{Y}_{s}\rightarrow \mathcal{X}_{s}$ is a good resolution. Write
    $$ 
    E_{s}= (p_{s}^{-1}\mathrm{Sing}(\mathcal{X}_{s}))_{\mathrm{red}} = \sum\limits_{i}E_{s,i}.
    $$
    \item There exists a distinguished component $Z$ of $\mathcal{Y}_{0}$ with a good resolution $p_{Z}: Z \rightarrow X$ whose reduced vertical boundary is
    $$
    B= \left( p_{Z}^{-1}\mathrm{Sing}(X) \right)_{\mathrm{red}}.
    $$
    Every other component of $\mathcal{Y}_{0}$ maps into $\mathrm{Sing}(X)$.
    \item After labeling the $E_{s,i}$, there are Cartier divisors $\mathcal{D}_{i}$ on $\mathcal{Y}$, obtained as total transforms of their horizontal closures, such that
    $$
    \mathcal{D}_{i}|_{\mathcal{Y}_{s}}=E_{s,i} \qquad(s\neq0).
    $$
    No $\mathcal{D}_{i}$ contains $Z$, and
    $$
    \mathrm{Supp}(\mathcal{D}_{i}|_{Z})\subseteq B.
    $$
    \end{enumerate}
\end{lemma}

\begin{proof}
    First, we choose a projective good resolution of the generic fibre which is a sequence of blow-ups along smooth centers. By \cite[proof of Lemma 4.44]{Kollar23Families}, after shrinking the punctured base and taking a finite base change, the good resolution on the generic fibre can be extended to a good resolution over the punctured base with labeled relative simple normal crossings exceptional divisor. Put $X_{sm}=X\setminus\mathrm{Sing}(X)$. Let $\mathcal{U} \subseteq \mathcal{X}$ be the locus where $\pi:\mathcal{X}\rightarrow T$ is smooth, and put $\Sigma=\mathcal{X}\setminus\mathcal{U}$. Then $\mathcal{U}\cap X=X_{sm}$. Choose the blow-up centers over the punctured base to lie over $\Sigma$. Extend the sequence by taking the scheme-theoretic closure of each center in the current total space and blowing up that closure. Every extended center remains in the inverse image of the closed subset $\Sigma$. The resulting modification is therefore an isomorphism over $\mathcal{U}$, and in particular over $X_{sm}$.

    Let $S=\mathrm{Sing}(X)$, regarded as a reduced closed subscheme of the central fibre in the total space, and let $\mathcal{I}_{S}$ be its ideal. On the extended modification, blow up the pullback of $\mathcal{I}_{S}$, and then resolve the total space together with the central fibre and the horizontal exceptional divisors. These operations may be chosen to preserve the given punctured family and the inverse image of $\mathcal{U}$.

    Apply semistable reduction to this model, retaining the horizontal divisors in the boundary and leaving the inverse image of $\mathcal{U}$ unchanged (cf. \cite[Theorem~2.58]{Kollar23Families}). Write the final reduced central fibre as
    $$
    \mathcal{Y}_{0}=Z+\sum\limits_jA_j,
    $$
    where $Z$ is the closure of the inverse image of $X_{sm}$. It is smooth, and $p_{Z}:Z\rightarrow X$ is projective, birational, and an isomorphism over $X_{sm}$. Every other component $A_j$ maps into $S$. The pullback of $\mathcal{I}_{S}$ remains invertible through the subsequent base changes and modifications. Hence
    $$
    \mathcal{I}_{S}\mathcal{O}_{\mathcal{Y}} = \mathcal{O}_{\mathcal{Y}}(-F)
    $$
    for an effective vertical Cartier divisor $F$. Here the ideal is pulled back through the composite morphism to the total space on which it was chosen. Since $Z$ dominates $X$ and every $A_j$ maps into $S$,
    $$
    \mathrm{Supp}(F)=\bigcup_jA_j.
    $$
    Taking the restriction on $Z$, we obtain
    $$
    \left(p_{Z}^{-1}(S)\right)_{\mathrm{red}} = \left(Z\cap\bigcup_jA_j\right)_{\mathrm{red}}.
    $$
    This is a simple normal crossings divisor on $Z$. Hence $p_{Z}$ is a good resolution with the asserted reduced vertical boundary.
    
    On the smooth model before semistable reduction, the horizontal exceptional divisors are Cartier. Retain their total pullbacks as the divisors $\mathcal{D}_{i}$. Their supports map into the relative non-smooth locus, which meets $X$ in $S$. Therefore no $\mathcal{D}_{i}$ contains $Z$, and
    $$
    \mathrm{Supp}(\mathcal{D}_{i}|_{Z})\subseteq B.
    $$
    Their restrictions to the noncentral fibres are unchanged.
\end{proof}

\begin{theorem}[Compatibility with $\mathbb{Q}$-smoothing]\label{thm:stationary-smoothing-compatibility}
    Notation for the family as in \cref{lem:extend-divisor-classes}, and for $g,n,\gamma_j,c_j,k_j$ as in \cref{def:stationary-resolution-invariant}. After shrinking the base, assume that $\mathcal{X}_{s}$ is terminal quotient Fano for $s \neq 0$, and denote its canonical covering stack by $\mathfrak{X}_{s}$. Let
    $$
    i_{0}:X\hookrightarrow\mathcal{X}, \qquad i_{s}: \mathcal{X}_{s} \hookrightarrow\mathcal{X}
    $$
    be the inclusions. Fix $\delta\in H_{2}(\mathcal{X},\mathbb{Z})$ and put $d=-K_{\mathcal{X}/T}\cdot\delta$. Assume $d\in\mathbb{Z}_{\geq 0}$ and
    $$
    \sum\limits_{j=1}^{n}(k_j+c_j)=d+n.
    $$
    When $d=0$, assume $2g-2+n>0$. Put $X_{sm}=X\setminus\mathrm{Sing}(X)$. Then, after shrinking $T$ and making an \'etale base change, there exist unique extensions $\gamma_{j,s}$ of $\gamma_{j}$ as untwisted insertions on $\mathfrak{X}_{s}$ such that, for sufficiently small $s \neq 0$, we have
    $$
    \sum\limits_{\substack{ \beta\in H_{2}(X,\mathbb{Z})\\ (i_{0})_{*}\beta=\delta}}
    \left\langle \prod_{j=1}^{n}\tau_{k_j}(\gamma_j) \right\rangle^{X,\mathrm{res}}_{g,n,\beta} =
    \sum\limits_{\substack{ \beta_{s}\in H_{2}(\mathcal{X}_{s},\mathbb{Z}) \\ (i_{s})_{*}\beta_{s}=\delta}}
    \left\langle \prod_{j=1}^{n}\tau_{k_j}(\gamma_{j,s}) \right\rangle^{\mathfrak{X}_{s}}_{g,n,\beta_{s}}.
    $$
    Both sums have finitely many nonzero terms.
\end{theorem}

\begin{proof}
    We first construct the extensions of the insertions. By \cref{lem:stationary-support}, $\gamma_{j}$ corresponds to a unique class in $H_{c}^{2c_{j}}(X_{sm},\mathbb{Q})$. By \cite[Corollary 6.15]{Siebenmann72}, over a sufficiently small analytic neighborhood of the central point, we can extend the class on the smooth locus and then extend by zero. For the finite collection of insertions, these choices can be made simultaneously. The resulting classes extend over the family, restrict to $\gamma_{j}$ on the central fibre, and are supported away from the central singular points. Since $\pi$ is projective, such an extension is unique in an analytic neighborhood of $0 \in T$. 

    We apply \cref{lem:scheme-resolution-family} to obtain a semistable reduction $h:\mathcal{Y}\rightarrow\mathcal{X}$ after the finite base change. Without loss of generality, we can still denote the new base by $T$. For $0 \neq s \in T$, write
    $$
    p_{s}:(\mathcal{Y}_{s},E_{s})\rightarrow \mathcal{X}_{s}, \qquad p:(Z,B)\rightarrow X
    $$
    for the general resolution and distinguished pair, with $E_{s}=\sum\limits_{i}E_{s,i}$ and $B=\sum\limits_{i}B_{i}$. Then
    $$
    h|_{\mathcal{Y}_{s}}=i_{s}\circ p_{s}, \qquad h|_{Z}=i_{0}\circ p.
    $$
    
    We compare resolution invariants, restricting to $\alpha_{s}\in H_{2}(\mathcal{Y}_{s},\mathbb{Z})$ satisfying
    $$
    (i_{s})_{*}(p_{s})_{*}\alpha_{s}=\delta, \qquad E_{s,i}\cdot\alpha_{s}=0 \quad\text{for every }i.
    $$
    Apply \cref{construction: deformation family 1} to $(\mathcal{Y}_{s},E_{s})$ with $A=p_{s}^{*}(-K_{\mathcal{X}_{s}})$. By \cref{lem:stationary-degeneration}, their logarithmic invariants equal the corresponding absolute invariants of $\mathcal{Y}_{s}$. Pull back $-K_{\mathcal{X}/T}$ to the semistable model. This gives a relatively nef divisor restricting to
    $$
    p_{s}^{*}(-K_{\mathcal{X}_{s}}) \quad\text{on }\mathcal{Y}_{s}, \qquad p^{*}(-K_{X}) \quad\text{on }Z.
    $$
    On classes with zero exceptional degree, $p_s^{*}(-K_{\mathcal{X}_s})$ and $-K_{\mathcal{Y}_s}$ have equal degrees. The extended stationary insertions vanish on auxiliary components and restrict to $p^{*}\gamma_{j}$ on $Z$. By \cref{lem:stationary-degeneration}, the nontrivial terms satisfy $B_{i}\cdot\alpha=0$ for every $i$ and have coefficient one. 
    
    Retain the horizontal exceptional divisor classes in the degeneration formula. Their restrictions to $Z$ are supported on $B$, so every surviving class satisfies the specialized generic exceptional-degree conditions. Retain the integral total-space classes and their pushforward $\delta$ under $h$ in the virtual identities. These conditions, together with the degrees of the extended exceptional divisors, are locally constant on the bounded mapping problems. Hence the virtual identities restrict to the corresponding open and closed substacks.

    The surviving central classes are precisely the integral $\alpha\in H_{2}(Z,\mathbb{Z})$ with
    $$
    B_{i}\cdot\alpha=0 \quad\text{for every }i, \qquad (i_0)_{*}p_{*}\alpha=\delta.
    $$
    For fixed $d$, \cref{thm:stationary-resolution-independence} gives finitely many contributing classes on the two resolutions. Apply the identities separately to their images in $H_{2}(\mathcal{Y},\mathbb{Z})$ and sum. A total-space class with no noncentral representative has zero complete central contribution by deformation invariance.
    
    Regrouping by $p_{*}\alpha$ and $(p_{s})_{*}\alpha_{s}$ gives the asserted identity with resolution invariants on both sides. Apply \cref{prop:quotient-resolution-comparison} to $\mathcal{X}_{s}$ to obtain its orbifold invariants.
\end{proof}

\subsection{Quantum periods of terminal Fano threefolds}\label{subsec:terminal-quantum-periods}

\begin{definition}[Quantum periods]\label{def:resolution-quantum-period}
    Let $X$ be a projective terminal Fano threefold, and let $\mathbf{1}\in H^{6}(X,\mathbb{Q})$ be the Poincar\'e dual of the point class. Define its quantum period and regularized quantum period by
    $$
    G_{X}^{\mathrm{res}}(t)
    = 1+\sum\limits_{d\geq2}t^d \sum\limits_{\substack{ \beta\in H_{2}(X,\mathbb{Z})\\
    -K_{X} \cdot\beta=d}} \left\langle \tau_{-K_{X} \cdot\beta-2}\mathbf{1} \right\rangle^{X,\mathrm{res}}_{0,1,\beta}
    $$
    and
    $$
    \widehat{G}_{X}^{\mathrm{res}}(t)
    = 1+\sum\limits_{d\geq2}d!\,t^d \sum\limits_{\substack{ \beta\in H_{2}(X,\mathbb{Z})\\
    (-K_{X})\cdot\beta=d}} \left\langle \tau_{-K_{X} \cdot\beta-2}\mathbf{1} \right\rangle^{X,\mathrm{res}}_{0,1,\beta}.
    $$
    
    Let $N^{1}(X)_{\mathbb{Q}}$ be the $\mathbb{Q}$-vector space generated by Cartier divisors modulo numerical equivalence, and set
    $$
    N^{1}(X)_{\mathbb{C}} = N^{1}(X)_{\mathbb{Q}} \otimes_{\mathbb{Q}}\mathbb{C}.
    $$
    For $D\in N^{1}(X)_{\mathbb{C}}$, define
    $$
    \widehat{G}_{X,D}^{\mathrm{res}}(t) = 1+\sum\limits_{d\geq2}d!\,t^d \sum\limits_{\substack{ \beta\in H_{2}(X,\mathbb{Z})\\ (-K_{X})\cdot\beta=d}} \left\langle \tau_{-K_{X} \cdot\beta-2}\mathbf{1} \right\rangle^{X,\mathrm{res}}_{0,1,\beta} e^{-D\cdot\beta}.
    $$
    By \cref{thm:stationary-resolution-independence}, the sum is finite for each $d$. The series are independent of the chosen resolution. We may omit the superscript $\mathrm{res}$ when there is no ambiguity.
\end{definition}

\begin{corollary}[Compatibility of quantum periods]\label[corollary]{cor:stationary-period-smoothing}
    Let $X$ be a projective terminal Fano threefold. For a projective $\mathbb{Q}$-smoothing with notation as in \cref{thm:stationary-smoothing-compatibility}, and sufficiently small $s\neq0$, we have
    $$
    G_{X}^{\mathrm{res}}(t)=G_{\mathcal{X}_{s}}(t),
    \qquad \widehat{G}_{X}^{\mathrm{res}}(t) = \widehat{G}_{\mathcal{X}_{s}}(t).
    $$
    Here $G_{\mathcal{X}_{s}}$ is the ordinary quantum period associated with the orbifold theory, obtained from its regularized period by dividing the degree-$d$ coefficient by $d!$.
    
    Let $D$ be a $\mathbb{Q}$-Cartier divisor on $X$, and let $m>0$ be an integer such that $mD$ is Cartier, and suppose that $\mathcal{O}_X(mD)$ extends to a line bundle $\mathcal{L}$ on $\mathcal{X}$. Put
    $$
    D_{s}= \frac{1}{m}c_{1}(\mathcal{L}|_{\mathcal{X}_{s}}).
    $$
    Then
    $$
    \widehat{G}_{X,D}^{\mathrm{res}}(t) = \widehat{G}_{\mathcal{X}_{s},D_{s}}(t).
    $$
    The same statement holds for complex linear combinations of extended $\mathbb{Q}$-Cartier divisor classes.
\end{corollary}

\begin{proof}
    Fix an integer $d\geq2$. Apply \cref{thm:stationary-smoothing-compatibility} with
    $$
    g=0,\qquad n=1,\qquad \gamma_1=\mathbf{1},\qquad k_1=d-2,
    $$
    and sum over the integral total-space classes $\delta$ with
    $$
    c_1(-K_{\mathcal{X}/T})\cdot\delta=d.
    $$
    Only finitely many contribute, and each fibre class occurs once, indexed by its total-space image. This gives the unweighted identity.
    
    For the weighted identity, put
    $$
    \theta=\frac{1}{m}c_1(\mathcal{L}) \in H^2(\mathcal{X},\mathbb{Q}).
    $$
    If $(i_0)_{*}\beta=(i_{s})_{*}\beta_{s}=\delta$, then
    $$
    D\cdot\beta=\theta\cdot\delta=D_{s}\cdot\beta_{s}.
    $$
    Thus the weights agree in each fixed-$\delta$ comparison. Summing and multiplying by $d!$ proves the regularized identities. The constant and degree-one terms agree by definition. The argument also applies to complex linear combinations of the extended divisor classes.
\end{proof}

\section{Quantum periods under the minimal model program} \label{sec:ordinary-terminal-case}

    In this section, we prove the main theorem. The following example shows that an individual flip need not preserve ordinary one-point Gromov-Witten invariants.

\begin{example}[A toric flip where the one-point invariant changes]\label{example: GW change after flip}
    Let
    $$
    R=\mathbb C[u,v,x_{0},x_{1},x_{2}]
    $$
    with \((\mathbb C^\ast)^2\)-weights
    $$
    \begin{array}{c|ccccc}
            & u & v & x_{0} & x_{1} & x_{2} \\ \hline
    F       & 1 & 1 & 0   & -1  & -2  \\
    H       & 0 & 0 & 1   & 1   & 1
    \end{array}
    $$
    Hence the torus invariant divisors are given by
    $$
    D_u=D_{v}=F,\qquad
    D_{x_{0}}=H,\qquad
    D_{x_{1}}=H-F,\qquad
    D_{x_{2}}=H-2F.
    $$
    In particular, the anticanonical divisor is given by
    $$
    -K=D_u+D_{v}+D_{x_{0}}+D_{x_{1}}+D_{x_{2}}=3H-F.
    $$
    
    Let $N\cong \mathbb Z^3$, and let the rays be
    $$
    \rho_{1}=(-1,1,2),\qquad
    \rho_{2}=(1,0,0),\qquad
    \rho_3=(0,-1,-1),
    $$
    $$
    \rho_4=(0,1,0),\qquad
    \rho_5=(0,0,1),
    $$
    corresponding respectively to
    $$
    u,\ v,\ x_{0},\ x_{1},\ x_{2} .
    $$
    They satisfy the two relations
    \[
    \rho_{1}+\rho_{2}-\rho_4-2\rho_5=0,
    \qquad
    \rho_3+\rho_4+\rho_5=0,
    \]
    which are dual to the two rows of the weight matrix above.
    
    We consider two adjacent GIT chambers:
    $$
    \sigma=\langle H,H-F\rangle,
    \qquad
    \sigma_{+}=\langle F,H\rangle.
    $$
    Let
    $$
    X:=X_{\sigma},
    \qquad
    X^{+}:=X_{\sigma_{+}}.
    $$
    The chamber $\sigma_{+}$ gives
    $$
    X^{+} = \mathbb{P}_{\mathbb{P}^{1}}
    \bigl(
    \mathcal{O}\oplus \mathcal{O}(1)\oplus \mathcal{O}(2)
    \bigr),
    $$
    with the quotient convention.
    
    We use the abbreviation
    $$
    \langle ijk\rangle:=\mathrm{Cone}(\rho_{i},\rho_{j},\rho_k).
    $$
    The maximal cones of the fans $\Sigma$ and $\Sigma_{+}$ are
    $$
    \Sigma(3)=
    \{
    \langle125\rangle,\,
    \langle235\rangle,\,
    \langle135\rangle,\,
    \langle124\rangle,\,
    \langle234\rangle,\,
    \langle134\rangle
    \}
    $$
    and
    $$
    \Sigma_{+}(3)=
    \{
    \langle245\rangle,\,
    \langle235\rangle,\,
    \langle234\rangle,\,
    \langle145\rangle,\,
    \langle135\rangle,\,
    \langle134\rangle
    \}.
    $$
    The varieties $X$ and $X^{+}$ are related by a toric flip. On the $X^{+}$-side, the flipped curve is
    $$
    C_{0}^{+}=V(\rho_4,\rho_5)=\{x_{1}=x_{2}=0\}\cong \mathbb P^{1}.
    $$
    Its normal bundle is
    $$
    N_{C_{0}^{+}/X^{+}}\cong \mathcal{O}_{\mathbb P^{1}}(-1)\oplus
    \mathcal{O}_{\mathbb P^{1}}(-2).
    $$
    On the \(X\)-side the flipping curve is
    $$
    C_{0}=V(\rho_{1},\rho_{2})=\{u=v=0\}.
    $$
    The only singular maximal cone of $\Sigma$ is
    $$
    \langle124\rangle,
    $$
    whose determinant is 2. Hence $X$ has one terminal quotient singularity of type $\frac{1}{2}(1,1,1)$. Moreover,
    $$
    \mathrm{Nef}(X)=\langle H,H-F\rangle,
    $$
    and
    $$
    -K_{X}=3H-F=2H+(H-F)
    $$
    lies in the interior of this cone. Hence $X$ is a terminal $\mathbb{Q}$-Fano threefold.
    
    Now let $L$ denote the class of a line in a fiber of
    $$
    \pi:X^{+} \rightarrow \mathbb{P}^{1}.
    $$
    Then
    $$
    F\cdot L=0,\qquad H\cdot L=1.
    $$
    The flipped curve $C_{0}^{+}$ satisfies
    $$
    F\cdot C_{0}^{+}=1,\qquad H\cdot C_{0}^{+}=0.
    $$
    Let $C^{+}\subset X^{+}$ be the section corresponding to the quotient
    $$
    \mathcal{O}\oplus \mathcal{O}(1)\oplus \mathcal{O}(2)
    \twoheadrightarrow
    \mathcal{O}(1).
    $$
    In toric notation,
    $$
    C^{+}=V(\rho_3,\rho_5)=\{x_{0}=x_{2}=0\}.
    $$
    This curve is disjoint from $C_{0}^{+}$. Its class is
    $$
    \beta^{+}:=[C^{+}]=[C_{0}^{+}]+[L],
    $$
    because both sides have the same intersections with $F$ and $H$:
    $$
    F\cdot \beta^{+}=1,\qquad H\cdot \beta^{+}=1.
    $$
    Therefore
    $$
    -K_{X^{+}}\cdot \beta^{+} = (3H-F)\cdot \beta^{+} = 2.
    $$
    The cone $\langle\rho_3,\rho_5\rangle$ belongs to both fans, so the strict transform of $C^{+}$ is
    $$
    C:=V(\rho_3,\rho_5)\subset X.
    $$
    Put $\beta=[C]\in H_{2}(X,\mathbb{Z})$. The two maximal cones of $\Sigma$ containing $\langle \rho_{3},\rho_{5}\rangle$ are $\langle 235 \rangle$ and $\langle 135 \rangle$. The corresponding wall relation is
    $$
    \rho_{3}-\rho_{5}+\rho_{1}+\rho_{2}=0.
    $$
    Hence
    $$
    D_{1}\cdot C=1,\qquad
    D_{2}\cdot C=1,\qquad
    D_3\cdot C=1,\qquad
    D_4\cdot C=0,\qquad
    D_5\cdot C=-1.
    $$
    Since
    $$
    D_{1}=D_{2}=F,\qquad D_3=H,\qquad D_4=H-F,\qquad D_5=H-2F,
    $$
    we get
    $$
    F\cdot C=1,\qquad H\cdot C=1,\qquad (H-F)\cdot C=0.
    $$
    By the definition of $\beta$, a stable map representing $\beta$ on $X$ is the normalization map
    $$
    f:\mathbb P^{1}\longrightarrow C=V(\rho_3,\rho_5)\subset X.
    $$
    This curve avoids the flipping locus $C_{0}=V(\rho_{1},\rho_{2})$.
    
    We now compute the one-point invariants in classes $\beta$ and $\beta^{+}$. Since
    $$
    -K_{X}\cdot\beta = -K_{X^{+}}\cdot\beta^{+} =2,
    $$
    we are computing
    $$
    \left\langle \tau_{0}\mathbf{1}\right\rangle_{\beta^{(+)}}
    = \left\langle \mathbf{1}\right\rangle_{\beta^{(+)}}.
    $$
    Choose a very general point
    $$
    p\in X^{+}\setminus D_{5}^{+},
    \qquad
    D_{5}^{+}=(x_{2}=0) \subset X^{+}.
    $$
    Let
    $$
    b=\pi(p),\qquad
    z_b=C_{0}^{+}\cap \pi^{-1}(b).
    $$
    The fiber
    $$
    \pi^{-1}(b)\cong \mathbb{P}^2
    $$
    contains a unique line $l_p$ through \(p\) and \(z_b\). Define a stable map $f_p:\Gamma_{1}\cup\Gamma_{2}\to X^{+}$
    by
    $$
    \Gamma_{1}\xrightarrow{\sim} C_{0}^{+},
    \qquad
    \Gamma_{2}\xrightarrow{\sim}l_p,
    $$
    with the node mapped to $z_b$, and with the marked point lying on $\Gamma_{2}$ and mapping to $p$. Then
    $$
    [f_p]=[C_{0}^{+}]+[l_p]=[C_{0}^{+}]+[L]=\beta^{+}.
    $$
    
    We claim that this is the unique stable map in class $\beta^{+}$ passing through $p$. Indeed, any stable map in class $\beta^{+}$ has
    $$
    F\cdot \beta^{+}=1,\qquad H\cdot \beta^{+}=1.
    $$
    If the component through $p$ were horizontal over $\mathbb{P}^{1}$, it would be a section of relative degree 1. Such a section corresponds to a quotient line bundle
    $$
    \mathcal{O} \oplus \mathcal{O}(1) \oplus \mathcal{O}(2)
    \twoheadrightarrow \mathcal{O}(1).
    $$
    However, we have $\mathrm{Hom}(\mathcal{O}(2),\mathcal{O}(1))=0$, so every such section lies in $D_{5}^{+}$. Since $p \notin D_{5}^{+}$, this is impossible. Hence the component through $p$ is vertical. Since $H\cdot\beta^{+}=1$, it is a line in the fiber $\pi^{-1}(b)$. The remaining horizontal class is then $C_{0}^{+}$, and the line is necessarily the unique line joining $p$ and $z_b$. Therefore
    $$
    ev^{-1}(p)=\{[f_p]\}.
    $$
    The virtual multiplicity of this point is 1. Indeed, we have
    $$
    T_{X^{+}}|_{C_{0}^{+}} \cong
    \mathcal{O}(2)\oplus \mathcal{O}(-1)\oplus \mathcal{O}(-2),
    $$
    $$
    T_{X^{+}}|_{l_p} \cong \mathcal{O}(2)\oplus \mathcal{O}(1)\oplus \mathcal{O}.
    $$
    The normalization sequence, together with the global generation of
    $T_{X^{+}}|_{l_p}$, gives
    $$
    H^{1}\!\left(
    \Gamma_{1}\cup\Gamma_{2},f_p^{*}T_{X^{+}}
    \right)
    \simeq
    H^{1}\!\left(
    C_{0}^{+},N_{C_{0}^{+}/X^{+}}
    \right)
    \simeq
    H^{1}\!\left(
    \mathbb{P}^{1},\mathcal{O}_{\mathbb{P}^{1}}(-2)
    \right)
    \simeq\mathbb{C}.
    $$
    This is the obstruction space for deformations of the map with the source fixed. We show that the node-smoothing map in the tangent-obstruction sequence is an isomorphism onto this space.
    
    Under the splitting
    $$
    N_{C_{0}^{+}/X^{+}}
    \simeq
    \mathcal{O}_{\mathbb{P}^{1}}(-1)
    \oplus
    \mathcal{O}_{\mathbb{P}^{1}}(-2),
    $$
    the first summand is tangent to $D_5^{+}$, while the second is normal to $D_5^{+}$. More precisely,
    $$
    N_{C_{0}^{+}/D_5^{+}}
    \simeq\mathcal{O}_{\mathbb{P}^{1}}(-1),
    \qquad
    N_{D_5^{+}/X^{+}}|_{C_{0}^{+}}
    \simeq\mathcal{O}_{\mathbb{P}^{1}}(-2).
    $$
    Since $p\notin D_5^{+}$, the line $l_p$ is not contained in the line $D_5^{+}\cap\pi^{-1}(b)$ in the fibre $\pi^{-1}(b)\simeq\mathbb{P}^{2}$. Hence $l_p$ meets $D_5^{+}$ transversely at $z_b$, and its tangent direction has nonzero projection to the $\mathcal{O}_{\mathbb{P}^{1}}(-2)$ normal direction.
    
    The node-smoothing map is the principal-part boundary map applied to this normal tangent direction. On the
    $\mathcal{O}_{\mathbb{P}^{1}}(-2)$ summand, it is induced by
    $$
    0\rightarrow
    \mathcal{O}_{\mathbb{P}^{1}}(-2)
    \rightarrow
    \mathcal{O}_{\mathbb{P}^{1}}(-2)(z_b)
    \rightarrow
    \left.
    \mathcal{O}_{\mathbb{P}^{1}}(-2)(z_b)
    \right|_{z_b}
    \rightarrow0.
    $$
    Indeed, locally a first-order smoothing of the node has equation $uv=\epsilon$. Substituting $v=\epsilon/u$ into the normal component of the map on the $l_p$-branch gives a simple pole on the $C_{0}^{+}$-branch, whose principal part is its normal tangent direction multiplied by $\epsilon$.
    
    Since
    $$
    \mathcal{O}_{\mathbb{P}^{1}}(-2)(z_b)
    \simeq\mathcal{O}_{\mathbb{P}^{1}}(-1)
    $$
    has vanishing $H^{0}$ and $H^{1}$, the connecting homomorphism
    $$
    H^{0}\!\left(
    z_b,
    \left.
    \mathcal{O}_{\mathbb{P}^{1}}(-2)(z_b)
    \right|_{z_b}
    \right)
    \xrightarrow{\sim}
    H^{1}\!\left(
    \mathbb{P}^{1},\mathcal{O}_{\mathbb{P}^{1}}(-2)
    \right)
    $$ is an isomorphism. The nonzero normal component found above therefore makes the node-smoothing map an isomorphism between one-dimensional spaces. Consequently, the obstruction space of the stable map vanishes, while the node-smoothing direction does not survive in its tangent space. The moduli space is smooth of the expected dimension at $[f_p]$. Moreover, as $q$ varies in $X^{+}\setminus D_{5}^{+}$, the same construction gives a family
    $$
    q\longmapsto [C_{0}^{+}\cup l_q,q],
    $$
    and the evaluation map is locally the identity on this family. Therefore \(ev\) has local degree 1 at
    $[f_p]$. Hence
    $$
    \left\langle \mathbf{1}\right\rangle^{X^+}_{\beta^+}=1.
    $$
    
    Now consider the flipping side \(X\). The corresponding class satisfies
    $$
    F\cdot\beta=H\cdot\beta=1.
    $$
    Therefore $(H-F)\cdot\beta=0$.
    Let
    $$
    D_5=(x_{2}=0),
    \qquad
    [D_5]=H-2F.
    $$
    For any non-constant component class \(\alpha\) of a stable map in total class $\beta$, nefness of $H-F$ gives
    \[
    (H-F)\cdot\alpha\ge 0.
    \]
    Since the total degree is zero,
    \[
    (H-F)\cdot\alpha=0
    \]
    for every non-constant component. Hence
    \[
    F\cdot\alpha=H\cdot\alpha.
    \]
    Since $H$ and $H-F$ span the nef cone and $2H+(H-F)=-K_{X}$ is ample, for every such class $\alpha$, we have $H \cdot \alpha>0$. Hence
    $$
    D_5\cdot\alpha = (H-2F)\cdot\alpha = -H\cdot\alpha < 0.
    $$
    An effective divisor has a negative intersection with a curve only if that curve is contained in the divisor. Hence every non-constant component of every stable map of class $\beta$ is contained in $D_{5}$.
    
    Choose the point insertion to be represented by a general smooth point $p\in X\setminus D_5$. Then no stable map of class $\beta$ can have its marked point mapping to \(p\). Hence $ev^{-1}(p)=\emptyset$ and $\left\langle \mathbf{1}\right\rangle^X_{\beta}=0$. Consequently
    $$
    \left\langle \tau_{0}\mathbf{1}\right\rangle^{X^{+}}_{\beta^{+}}=1,
    \qquad
    \left\langle \tau_{0}\mathbf{1}\right\rangle^{X}_{\beta}=0.
    $$
    Hence the one-point descendant Gromov-Witten invariant is not invariant under this flip.
\end{example}

\subsection{Generalized relative Gromov-Witten invariants}\label{subsec:ordinary-resolution-independence}

    We specify the geometric data for the auxiliary relative comparison.

\begin{statement}[Resolutions used for the relative comparison]\label{def:ordinary-admissible}
    Let $V$ be a projective $\mathbb{Q}$-factorial terminal threefold, and let $D$ be an effective $\mathbb{Q}$-Cartier divisor, possibly zero. We consider a projective birational morphism $p: Z\rightarrow V$ from a smooth projective threefold, together with a reduced simple normal crossings divisor $B=\sum\limits_{i}B_{i}$. The divisor $B$ contains every exceptional prime divisor of $p$ and the support of $p^{*}D$. Every $B_{i}$ is exceptional over $V$ or maps onto a component of $D$. Write
    $$
    K_{Z}=p^{*}K_{V}+\sum\limits_{i} a_{i}B_{i}, \qquad p^{*}D=\sum\limits_{i} m_{i}B_{i}.
    $$
    Here $m_{i}\geq 0$, and
    $$
    a_{i}>0 \quad\text{for exceptional }B_{i}, \qquad a_{i}=0 \quad\text{for nonexceptional }B_{i}.
    $$
    We permit further blow-ups along smooth centers of codimension at least two having normal crossings with the boundary.
\end{statement}

    In this subsection, we generalize the descendant point-insertion relative Gromov-Witten invariants to terminal Fano threefolds.

\begin{definition}[Relative Gromov-Witten invariants]\label{def:ordinary-resolved-coefficient}
    Notation as in \cref{def:ordinary-admissible}. Fix a rational number $c\geq 0$ such that the divisor
    $$
    A_{V} := -K_{V}+cD
    $$
    is ample, with either $c>0$ or $D=0$. Choose a general point $q\in Z\setminus B$. For $\beta\in H_{2}(V,\mathbb{Z})$ satisfying
    $$
    D\cdot\beta=0, \qquad d=(-K_{V})\cdot\beta\in\mathbb{Z}_{\geq2},
    $$
    we define
    $$
    \left\langle \tau_{d-2}\mathbf{1}\mid\emptyset \right\rangle^{(V,D),p}_{0,1,\beta} = \sum\limits_{\substack{
    \alpha\in H_{2}(Z,\mathbb{Z})\\
    p_{*}\alpha=\beta\\
    B_{i}\cdot\alpha=0\ \text{for every }i}}
    \left\langle \tau_{d-2}\mathbf{1}\mid\emptyset \right\rangle^{(Z,B),\log}_{0,1,\alpha}.
    $$
    The relative moduli space has one original marking with contact vector zero and no relative markings.
\end{definition}

\begin{lemma}[Finiteness]\label{lem:ordinary-finiteness}
    Notation as in \cref{def:ordinary-resolved-coefficient}. Then the sum defining $\left\langle \tau_{d-2}\mathbf{1}\mid\emptyset \right\rangle^{(V,D),p}_{0,1,\beta}$ has finitely many nonzero terms. Moreover, for fixed $p$, $c$, and $d$, only finitely many integral classes $\alpha\in H_{2}(Z,\mathbb{Z})$ contribute as $\beta$ varies subject to
    $$
    A_{V}\cdot\beta=d, \qquad D\cdot\beta=0.
    $$
\end{lemma}

\begin{proof}
    Choose an ample $\mathbb{Q}$-divisor $H$ on $Z$. Since $V$ is $\mathbb{Q}$-factorial, $H_{V}=p_{*}H$ is $\mathbb{Q}$-Cartier, and we can write
    $$
    H = p^{*}H_{V} + F
    $$
    where $F$ is supported on exceptional divisors. In particular, we have $F \cdot \alpha = 0$. Hence, by the projection formula, we have
    $$
    H \cdot \alpha = H_{V} \cdot \beta.
    $$
    Such curve classes $\alpha$ are finite. For the second assertion, choose $a>0$ such that $aA_{V} - H_{V}$ is ample. Then
    $$
    H\cdot\alpha=H_{V}\cdot\beta\leq ad.
    $$
    After replacing $H$ by a Cartier multiple, these are bounds on an ample Cartier degree. The relevant moduli spaces of stable maps are of finite type, and the integral homology class is locally constant. Hence only finitely many integral classes occur.
\end{proof}

\begin{proposition}[Resolution independence]\label[proposition]{prop:ordinary-resolution-independence}
    Notation as in \cref{def:ordinary-resolved-coefficient}. Then for fixed $(V,D)$ and $\beta$, the value $\left\langle \tau_{d-2}\mathbf{1}\mid\emptyset \right\rangle^{(V,D),p}_{0,1,\beta}$ is independent of the resolution $p$. In particular, we can denote this common value by $\left\langle \tau_{d-2}\mathbf{1}\mid\emptyset \right\rangle^{(V,D)}_{0,1,\beta}$.
\end{proposition}

\begin{proof}
    The proof is similar to the proof of \cref{thm:stationary-resolution-independence}. We first consider a single blow-up $f: \widetilde{Z} \rightarrow Z$ along a smooth center of codimension at least two having normal crossings with $B$. Fix a curve class $\alpha$ such that $p_{*}\alpha=\beta$ and $B_{i}\cdot\alpha=0$ for every $i$. Then
    $$
    (-K_{Z})\cdot\alpha = (-K_{V})\cdot\beta = A_{V}\cdot\beta = d.
    $$
    
    Apply \cref{construction: deformation family 1} to $(Z,B)$ with $A=p^{*}A_{V}$. By \cref{lem:ordinary-marked-factor,lem:ordinary-edge-free}, we have
    $$
    \left\langle \tau_{d-2}\mathbf{1} \right\rangle^Z_{0,1,\alpha} =
    \left\langle \tau_{d-2}\mathbf{1}\mid\emptyset \right\rangle^{(Z,B),\log}_{0,1,\alpha}.
    $$
    Apply \cref{construction: deformation family 2} to $f: \widetilde{Z} \rightarrow Z$ with $A=p^{*}A_{V}$. Write $\widetilde{B}=\sum\limits_{j}\widetilde{B}_{j}$ for the reduced union of the new exceptional divisor and the strict transforms of the components of $B$. The same vanishing argument gives
    $$
    \left\langle \tau_{d-2}\mathbf{1} \right\rangle^Z_{0,1,\alpha} =
    \sum\limits_{\substack{ f_{*}\widetilde{\alpha}=\alpha\\
    \widetilde{B}_{j}\cdot\widetilde{\alpha}=0\ \forall j}}
    \left\langle \tau_{d-2}\mathbf{1}\mid\emptyset \right\rangle^{(\widetilde{Z},\widetilde{B}),\log}_{0,1,\widetilde{\alpha}}.
    $$
    Every $f^{*}B_{i}$ is supported on $\widetilde{B}$, so the new boundary equalities imply the old ones. The morphisms of both total families to $Z$ retain $\alpha\in H_{2}(Z,\mathbb{Z})$ itself. As in \cref{lem:stationary-blowup}, restrict the virtual identities to the open and closed substacks with this integral pushforward. The surviving classes satisfy $f_{*}\widetilde{\alpha}=\alpha$ in integral homology. Sum over the finitely many relevant total-space homology classes. A class with no generic representative has zero complete special contribution by deformation invariance. Finiteness follows from \cref{lem:ordinary-finiteness}.
    
    For two resolutions, we apply the weak factorization as in the proof of \cref{thm:stationary-resolution-independence}, and the result follows.
\end{proof}

\begin{remark}
    One can define the relative point-insertion Gromov-Witten invariants for $n$ marked points and arbitrary genus $g$ in a similar way. However, resolution independence fails for the general insertion class setting as in \cref{sec:stationary-terminal-invariants}.
\end{remark}

\subsection{Resolution coefficients and nonpositive divisor degree}\label{subsec:terminal-resolution-coefficients}

\begin{proposition}[Identification of the one-point coefficient]\label[proposition]{prop:ordinary-absolute-resolved}
    Let $V$ be a projective $\mathbb{Q}$-factorial terminal Fano threefold, and let $\beta\in H_{2}(V,\mathbb{Z})$ satisfy
    $$
    d=(-K_V)\cdot\beta\in\mathbb{Z}_{\geq2}.
    $$
    Then
    $$
    \left\langle \tau_{d-2}\mathbf{1} \right\rangle^{V}_{0,1,\beta}=\left\langle \tau_{d-2}\mathbf{1}\mid\emptyset \right\rangle^{(V,0)}_{0,1,\beta}.
    $$
\end{proposition}

\begin{proof}
    On a good resolution with its reduced exceptional boundary, the two defining sums coincide. \cref{prop:ordinary-resolution-independence} permits computation on any resolution used in the relative comparison, including further modifications over the smooth locus.
\end{proof}

\begin{proposition}[Nonpositive divisor degree]\label[proposition]{prop:ordinary-nonpositive-degree}
    Let $V$, $\beta$, and $d$ be as in \cref{prop:ordinary-absolute-resolved}. Let $D$ be an effective $\mathbb{Q}$-Cartier divisor on $V$. Then
    $$
    \left\langle \tau_{d-2}\mathbf{1} \right\rangle^{V}_{0,1,\beta}=0 \quad\text{if }D\cdot\beta<0,
    $$
    and
    $$
    \left\langle \tau_{d-2}\mathbf{1} \right\rangle^{V}_{0,1,\beta}=\left\langle \tau_{d-2}\mathbf{1}\mid\emptyset \right\rangle^{(V,D)}_{0,1,\beta} \quad\text{if }D\cdot\beta=0.
    $$
\end{proposition}

\begin{proof}
    Choose $c \in \mathbb{Q}_{>0}$ sufficiently small such that
    $$
    L_{V}=-K_{V}+cD
    $$
    is ample. Take a projective log resolution $p:Z\rightarrow V$ of the singular locus and the support of $D$, such that
    $$
    F=\mathrm{Exc}(p)_{\mathrm{red}}=\sum\limits_{i}F_{i}, \qquad B= \left( F+\mathrm{Supp}(p^{*}D) \right)_{\mathrm{red}} = \sum\limits_{j}B_{j}
    $$
    are simple normal crossing divisors, and choose a general point $q \in Z\setminus B$. By the preceding proposition and resolution independence,
    $$
    \left\langle \tau_{d-2}\mathbf{1} \right\rangle^{V}_{0,1,\beta} = \sum\limits_{\substack{ p_{*}\alpha=\beta\\
    F_{i}\cdot\alpha=0\ \forall i}} \left\langle \tau_{d-2}\mathbf{1}\mid\emptyset \right\rangle^{(Z,F),\log}_{0,1,\alpha}.
    $$
    For each class, apply \cref{construction: deformation family 1} to $(Z,F)$ with $A=p^{*}(-K_{V})$. By \cref{lem:ordinary-marked-factor,lem:ordinary-edge-free}, this identifies the relative invariant with
    $$
    \left\langle \tau_{d-2}\mathbf{1} \right\rangle^Z_{0,1,\alpha}.
    $$
    
    Apply \cref{construction: deformation family 1} again, now to $(Z,B)$ with $A=p^{*}L_{V}$. The degree of $\alpha$ against $p^{*}L_{V}$ is
    $$
    d+ce, \qquad e=D\cdot\beta.
    $$
    \cref{lem:ordinary-marked-factor} gives zero for $e<0$. For $e=0$, only the contribution with no relative marked points remains.
    
    By \cref{lem:ordinary-edge-free}, that contribution is the logarithmic invariant of $(Z,B)$ in class $\alpha$ if $B_{j}\cdot\alpha=0$ for every $j$, and is zero otherwise.

    The class is exactly $\alpha$, because the total family maps to $Z$ and restricts to the identity on the distinguished component.
    
    Summing proves the assertion. The sums are finite by \cref{lem:ordinary-finiteness}.
\end{proof}

\subsection{The main theorem}\label{subsec:ordinary-period-limit}

    We now prove the main theorem.

\begin{theorem}\label{thm:ordinary-period-limit}
    Let $X$ and $X'$ be projective $\mathbb{Q}$-factorial terminal Fano threefolds. Let
    $$
    g:X\dashrightarrow X'
    $$
    be a birational map which contracts exactly one prime divisor $E$ and extracts no prime divisor. Then
    $$
    \lim_{r\rightarrow+\infty} \widehat{G}_{X,rE}^{\mathrm{res}}(t) = \widehat{G}_{X'}^{\mathrm{res}}(t).
    $$
    Here $\widehat{G}_{X,rE}^{\mathrm{res}}(t)$ is the regularized quantum period of $(X,rE)$ and $\widehat{G}_{X'}^{\mathrm{res}}(t)$ is the regularized quantum period of $X'$.
\end{theorem}

\begin{proof}
    Choose a common smooth projective resolution
    $$
    p:Z\rightarrow X, \qquad q:Z\rightarrow X',
    $$
    such that $q=g\circ p$ as rational maps and the union of the exceptional divisors and the strict transform $\widetilde{E}$ has simple normal crossings.

    Since $g$ extracts no prime divisor, every $p$-exceptional prime divisor is $q$-exceptional. Conversely, a $q$-exceptional prime divisor which is not $p$-exceptional is the strict transform of a prime divisor contracted by $g$, and hence equals $\widetilde{E}$. Thus the reduced divisor
    $$
    B:= \widetilde{E}+\mathrm{Exc}(p)_{\mathrm{red}} = \mathrm{Exc}(q)_{\mathrm{red}}
    $$
    has simple normal crossings. Write $B=\sum\limits_{i=1}^{k}B_{i}$. We use the same pair $(Z,B)$ to compute $\left\langle \tau_{d-2}\mathbf{1}\mid\emptyset \right\rangle^{(X,E)}_{0,1,\beta}$ and $\left\langle \tau_{d-2}\mathbf{1}\mid\emptyset \right\rangle^{(X',0)}_{0,1,\beta'}$ by \cref{prop:ordinary-resolution-independence}. For $(X,E)$, choose $c\in\mathbb{Q}_{>0}$ sufficiently small that $-K_{X}+cE$ is ample. For $(X',0)$, use the ample divisor $-K_{X'}$. Choose the point insertion in a common open subset where both $p$ and $q$ are isomorphisms. The divisors
    $$
    K_{Z}-p^{*}K_{X}, \qquad K_{Z}-q^{*}K_{X'}, \qquad p^{*}E
    $$
    are supported on $B$. Consequently, for every $\alpha\in H_2(Z,\mathbb{Z})$ satisfying $B_{i}\cdot\alpha=0$ for all $i$, we have $E\cdot p_{*}\alpha=0$ and
    $$
    (-K_{X})\cdot p_{*}\alpha = (-K_{Z})\cdot\alpha = (-K_{X'})\cdot q_{*}\alpha.
    $$
    Fix an integer $d\geq2$. By \cref{prop:ordinary-nonpositive-degree},
    $$
    \sum\limits_{\substack{ \beta\in H_{2}(X,\mathbb{Z})\\ (-K_{X})\cdot\beta=d\\ E\cdot\beta=0}} \left\langle \tau_{d-2}\mathbf{1} \right\rangle^{X}_{0,1,\beta} =
    \sum\limits_{\substack{ \beta\in H_{2}(X,\mathbb{Z})\\ (-K_{X})\cdot\beta=d\\ E\cdot\beta=0}} \left\langle \tau_{d-2}\mathbf{1}\mid\emptyset \right\rangle^{(X,E)}_{0,1,\beta}.
    $$
    By \cref{prop:ordinary-resolution-independence}, the latter sum can be computed on $(Z,B)$:
    $$
    \sum\limits_{\substack{\beta\in H_{2}(X,\mathbb{Z})\\ (-K_{X})\cdot\beta=d\\ E\cdot\beta=0}} \left\langle \tau_{d-2}\mathbf{1}\mid\emptyset \right\rangle^{(X,E)}_{0,1,\beta} =
    \sum\limits_{\substack{ \alpha\in H_{2}(Z,\mathbb{Z})\\ (-K_{Z})\cdot\alpha=d\\ B_{i}\cdot\alpha=0\ (1\leq i\leq k)}}
    \left\langle \tau_{d-2}\mathbf{1}\mid\emptyset \right\rangle^{(Z,B),\log}_{0,1,\alpha}.
    $$
    Computing $\left\langle \tau_{d-2}\mathbf{1}\mid\emptyset \right\rangle^{(X',0)}_{0,1,\beta'}$ on the same pair and applying \cref{prop:ordinary-absolute-resolved} gives
    $$
    \begin{aligned}
    \sum\limits_{\substack{\alpha\in H_{2}(Z,\mathbb{Z})\\ (-K_{Z})\cdot\alpha=d\\ B_{i}\cdot\alpha=0\ (1\leq i\leq k)}}
    \left\langle \tau_{d-2}\mathbf{1}\mid\emptyset \right\rangle^{(Z,B),\log}_{0,1,\alpha} &=
    \sum\limits_{\substack{\beta'\in H_{2}(X',\mathbb{Z})\\ (-K_{X'})\cdot\beta'=d}} \left\langle \tau_{d-2}\mathbf{1}\mid\emptyset \right\rangle^{(X',0)}_{0,1,\beta'}\\
    &= \sum\limits_{\substack{\beta'\in H_{2}(X',\mathbb{Z})\\ (-K_{X'})\cdot\beta'=d}} \left\langle \tau_{d-2}\mathbf{1} \right\rangle^{X'}_{0,1,\beta'}.
    \end{aligned}
    $$
    The sums are finite by \cref{lem:ordinary-finiteness} and regroup the same classes on $Z$ by their pushforwards $p_*\alpha$ and $q_*\alpha$.

    \cref{prop:ordinary-nonpositive-degree} also gives $\left\langle \tau_{d-2}\mathbf{1} \right\rangle^{X}_{0,1,\beta}=0$ whenever $E\cdot\beta<0$. There are only finitely many nonzero coefficients of anticanonical degree $d$, by \cref{thm:stationary-resolution-independence}. Thus the weights $e^{-rE\cdot\beta}$ remove precisely the positive-$E$-degree terms as $r\rightarrow+\infty$. Consequently,
    $$
    \begin{aligned}
    \lim_{r\rightarrow+\infty}
    [t^d]\widehat{G}_{X,rE}^{\mathrm{res}}(t)
    &=
    d!\sum\limits_{\substack{
    (-K_{X})\cdot\beta=d\\
    E\cdot\beta=0}}
    \left\langle \tau_{d-2}\mathbf{1} \right\rangle^{X}_{0,1,\beta}\\
    &=
    d!\sum\limits_{(-K_{X'})\cdot\beta'=d}
    \left\langle \tau_{d-2}\mathbf{1} \right\rangle^{X'}_{0,1,\beta'}\\
    &=
    [t^d]\widehat{G}_{X'}^{\mathrm{res}}(t).
    \end{aligned}
    $$
    Both constant terms equal one, and both coefficients of $t$ are zero by definition. This proves the theorem.
\end{proof}

\printbibliography

\end{document}